\documentclass[dvipsnames]{article}

\usepackage{xcolor}
\definecolor{darkred}{rgb}{0.45,0,0}
\usepackage[pagebackref,colorlinks,linkcolor=darkred,citecolor=darkred]{hyperref}

\usepackage{amsmath,amssymb,latexsym,mathpartir,fontawesome,amsthm,cleveref,url,
rotating,multirow,stmaryrd,eulervm,charter,xcolor,mathtools,graphicx,stackengine,pifont,tikz-cd,bbold,bbm,scalerel}

\usepackage[minimal]{leadsheets}
\useleadsheetslibraries{musicsymbols}

\usepackage{fancybox}

\usepackage{geometry}
\usepackage[title]{appendix}

\usepackage{mdframed}

\newmdenv[
  topline=false,
  bottomline=false,
  rightline=false,
  skipabove=\topsep,
  skipbelow=\topsep
]{siderules}

\usepackage{sectsty}
\sectionfont{\centering \large}
\subsectionfont{\centering \normalfont \textsc}

\usetikzlibrary{decorations.pathreplacing, decorations.markings, shapes, decorations.pathmorphing, patterns, arrows}

\tikzset{
  on each segment/.style={
    decorate,
    decoration={
      show path construction,
      moveto code={},
      lineto code={
        \path [#1]
        (\tikzinputsegmentfirst) -- (\tikzinputsegmentlast);
      },
      curveto code={
        \path [#1] (\tikzinputsegmentfirst)
        .. controls
        (\tikzinputsegmentsupporta) and (\tikzinputsegmentsupportb)
        ..
        (\tikzinputsegmentlast);
      },
      closepath code={
        \path [#1]
        (\tikzinputsegmentfirst) -- (\tikzinputsegmentlast);
      },
    },
  },
  mid arrow/.style={postaction={decorate,decoration={
        markings,
        mark=at position .5 with {\arrow[scale=1.5, #1]{to}}
      }}},
}

\newcommand*\circled[1]{\tikz[baseline=(char.base)]{
            \node[shape=circle,draw,thick,minimum size=2em,inner sep=2pt] (char) {#1};}}

\newcommand{\slfrac}[2]{\left.#1\middle/#2\right.}

\def\n{\cblu{n}}
\def\k{\cblu{k}}
\def\l{\cblu{l}}
\def\m{\cblu{m}}
\def\p{\cblu{p}}
\def\a{\cblu{a}}
\def\b{\cblu{b}}
\def\rc{\cblu{c}}
\def\A{\cblu{A}}
\def\B{\cblu{B}}
\def\C{\cblu{C}}
\def\D{\cblu{D}}
\def\E{\cblu{E}}
\def\F{\cblu{F}}
\def\G{\cblu{G}}
\def\x{\cblu{x}}
\def\y{\cblu{y}}
\def\z{\cblu{z}}
\def\X{\cblu{X}}
\def\Y{\cblu{Y}}
\def\c{\cblu{\mu}}
\def\i{\cblu{i}}
\def\j{\cblu{j}}
\def\k{\cblu{k}}
\def\SST{\mathsf{\cred{SST}}}
\def\Z{\cred{\mathsf{Z}}}
\def\S{\cred{\mathsf{S}}}
\def\Type{\mathsf{\cred{Type}}}
\def\Top{\mathsf{\cred{Top}}}
\def\d{{^\mathsf{\cred{d}}}}
\def\DD{{^\mathsf{\cred{D}}}}
\def\Path{\mathsf{\cred{Path}}}
\def\PathP{\mathsf{\cred{PathP}}}
\def\Perm{\mathsf{\cred{Perm}}}
\def\t{\cblu{t}}
\def\s{\cblu{s}}
\def\f{\cblu{f}}
\def\g{\cblu{g}}
\def\a{\cblu{a}}
\def\b{\cblu{b}}
\def\Base{\mathsf{\cred{Base}}}
\def\sq{{\cred{\square}}}
\def\lock{\text{\faLock}}
\def\ctx{\;\mathsf{ctx}}
\def\tel{\;\mathsf{tel}}
\def\term{\;\mathsf{term}}
\def\ev{{^\mathsf{\cred{ev}}}}
\def\I{\cred{\mathbb{I}}}
\def\R{\cred{\mathsf{R}}}
\def\N{\cred{\mathbb{N}}}

\let\oldlambda\lambda
\renewcommand{\lambda}{\cred{\oldlambda}}

\let\oldGamma\Gamma
\renewcommand{\Gamma}{\cblu{\oldGamma}}

\let\oldgamma\gamma
\renewcommand{\gamma}{\cblu{\oldgamma}}

\let\oldDelta\Delta
\renewcommand{\Delta}{\cblu{\oldDelta}}

\let\oldsigma\sigma
\renewcommand{\sigma}{\cblu{\oldsigma}}

\let\oldalpha\alpha
\renewcommand{\alpha}{\cblu{\oldalpha}}

\let\oldbeta\beta
\renewcommand{\beta}{\cblu{\oldbeta}}

\let\olddelta\delta
\renewcommand{\delta}{\cblu{\olddelta}}

\def\pair#1#2{\cred{\boldsymbol\langle}\; #1\;\cred{\boldsymbol,}\;#2 \;\cred{\boldsymbol\rangle}}

\definecolor{codegreen}{rgb}{0,0.6,0}
\definecolor{codegray}{rgb}{0.5,0.5,0.5}
\definecolor{codepurple}{rgb}{0.58,0,0.84}
\definecolor{backcolour}{rgb}{0.95,0.95,0.92}

\usepackage{listings} 
\makeatletter
\newcommand{\morphic@sqcdot}[2]{%
  \sbox\z@{$\m@th#1\centerdot$}%
  \ht\z@=.33333\ht\z@
  \vcenter{\box\z@}%
}
\renewcommand{\dot}{\mathbin{\mathpalette\morphic@sqcdot\relax}}
\makeatother

\title{You Wouldn't Permutahedron \vspace{-0.2cm}}

\author{Astra Kolomatskaia}

\date{\vspace{-3ex}}

\definecolor{darkgreen}{rgb}{0,0.4,0}
\newcommand{\cblu}[1]{{\color{blue}#1}}
\newcommand{\cred}[1]{{\color{darkgreen}#1}}
\newcommand{\cprime}[1]{{\color{red}#1}}

\makeatletter
\newsavebox{\@brx}
\newcommand{\llangle}[1][]{\savebox{\@brx}{\(\m@th{#1\cred{\boldsymbol\langle}}\)}%
  \mathopen{\copy\@brx\kern-0.5\wd\@brx\usebox{\@brx}}}
\newcommand{\rrangle}[1][]{\savebox{\@brx}{\(\m@th{#1\cred{\boldsymbol\rangle}}\)}%
  \mathclose{\copy\@brx\kern-0.5\wd\@brx\usebox{\@brx}}}
\makeatother

\newcommand{\dangle}[1]{\llangle~#1~\rrangle}

\makeatletter
\newsavebox{\@bry}
\newcommand{\newllparenthesis}[1][]{\savebox{\@bry}{\(\m@th{#1\cred{\boldsymbol(}}\)}%
  \mathopen{\copy\@bry\kern-0.5\wd\@bry\usebox{\@bry}}}
\newcommand{\newrrparenthesis}[1][]{\savebox{\@bry}{\(\m@th{#1\cred{\boldsymbol)}}\)}%
  \mathclose{\copy\@bry\kern-0.5\wd\@bry\usebox{\@bry}}}
\makeatother

\makeatletter
\newsavebox{\@brz}
\newcommand{\newllbracket}[1][]{\savebox{\@brz}{\(\m@th{#1\cred{\boldsymbol[}}\)}%
  \mathopen{\copy\@brz\kern-0.5\wd\@brz\usebox{\@brz}}}
\newcommand{\newrrbracket}[1][]{\savebox{\@brz}{\(\m@th{#1\cred{\boldsymbol]}}\)}%
  \mathclose{\copy\@brz\kern-0.5\wd\@brz\usebox{\@brz}}}
\makeatother

\newcommand{\dsqbkt}[1]{\newllbracket#1\newrrbracket}

\newcommand{\blank}{\mathord{\hspace{1pt}\text{--}\hspace{1pt}}}

\newcommand{\set}[1]{\left\{#1\right\}}
\newcommand{\bkt}[1]{\left(#1\right)}
\newcommand{\sqbkt}[1]{[#1]}

\newcommand{\sub}[1]{_{\scalebox{0.55}[1.0]{$\scriptscriptstyle #1$}}}
\newcommand{\ff}{\mathfrak{f}}


\newcommand{\bbox}{\rule{0em}{1ex}\hfill\ensuremath{\triangleleft}}

\newsavebox{\mybox}

\begin{document}

\maketitle

\begin{abstract}
\noindent
We develop formulas that define permutahedral commutation coherence relations of all orders. To illustrate the result geometrically, we begin by defining a rigid transformation of the $(\n+\cred{1})$-permutahedron into a $\n$-cube of dimensions $\cred{1} \times \cred{2} \times \cdots \times \n$. With a fictitious assumption, we then \emph{`define'} the corresponding coherence relations \emph{`up to associativity'} as an instance of a semi-simplicial type in the language of \emph{Displayed Type Theory}. This is not a formal result in type theory, but we expect this to translate into one as soon as the problem of defining associahedral coherences is solved in a type theory with semi-simplicial types. On the other hand, this not-strictly-well-typed definition may be used to produce well-typed formulas in a restricted setting.
\end{abstract}

\vspace{0.3cm}

\noindent
\emph{`Somebody is eventually going to write down a coherence for the \cred{$4$}-permutahedron, but I can say with confidence that nobody will ever write down a coherence for the \cred{$5$}-permutahedron.'}

\noindent
-- \textbf{the author}, \emph{circa ten months ago}

\vspace{0.2cm}

\begin{center}
{\cblu{\circled{\footnotesize \cred{\trebleclef}}}\ \  {\large \textbf{Prelude}}}
\end{center}

The permutahedral coherence relations share a close connection with the problem of constructing semi-simplicial types (SSTs), and we presently motivate them as such. A semi-simplicial type is intuitively a semi-simplicial set valued in homotopy types, and there are two non-rigorous notions specifying what this should mean.

The first, known as the \emph{indexed} formulation, says that a $\SST$ consists of an infinite list of the data that starts out as follows:
\begin{align*}
  \cred{A_0} :~&\Type \\
  \cred{A_1} :~&{\normalcolor (\cblu{x\sub{01}} :
  \cred{A_0})}~{\normalcolor (\cblu{x\sub{10}} :
  \cred{A_0})} \to \Type \\
  \cred{A_2} :~&{\normalcolor (\cblu{x\sub{001}} :
  \cred{A_0})}~{\normalcolor (\cblu{x\sub{010}} :
  \cred{A_0})}
  ~{(\cblu{\beta\sub{011}} : \cred{A_1}~\cblu{x\sub{001}~x\sub{010}})}
  ~{\normalcolor (\cblu{x\sub{100}} : \cred{A_0})}
  ~{(\cblu{\beta\sub{101}} : \cred{A_1}~\cblu{x\sub{001}~x\sub{100}})}
  ~{(\cblu{\beta\sub{110}} : \cred{A_1}~\cblu{x\sub{010}~x\sub{100}})} \to
  \Type \\
  \cred{A_3} :~&{\normalcolor (\cblu{x\sub{0001}} :
  \cred{A_0})}~{\normalcolor (\cblu{x\sub{0010}} :
  \cred{A_0})}
  ~{(\cblu{\beta\sub{0011}} : \cred{A_1}~\cblu{x\sub{0001}~x\sub{0010}})}
  ~{\normalcolor (\cblu{x\sub{0100}} : \cred{A_0})}
  ~{(\cblu{\beta\sub{0101}} : \cred{A_1}~\cblu{x\sub{0001}~x\sub{0100}})}
  ~{(\cblu{\beta\sub{0110}} : \cred{A_1}~\cblu{x\sub{0010}~x\sub{0100}})} \\
  &{(\cblu{\ff\sub{0111}} :
  \cred{A_2}~\cblu{x\sub{0001}~x\sub{0010}~\beta\sub{0011}~x\sub{0100}~\beta\sub{0101}~\beta\sub{0110}})}
  ~{\normalcolor (\cblu{x\sub{1000}} : \cred{A_0})}
  ~{(\cblu{\beta\sub{1001}} : \cred{A_1}~\cblu{x\sub{0001}~x\sub{1000}})}
  ~{(\cblu{\beta\sub{1010}} : \cred{A_1}~\cblu{x\sub{0010}~x\sub{1000}})}\\
  &{(\cblu{\ff\sub{1011}} :
  \cred{A_2}~\cblu{x\sub{0001}~x\sub{0010}~\beta\sub{0011}~x\sub{1000}~\beta\sub{1001}~\beta\sub{1010}})}
  ~{(\cblu{\beta\sub{1100}} : \cred{A_1}~\cblu{x\sub{0100}~x\sub{1000}})}
  ~{(\cblu{\ff\sub{1101}} :
  \cred{A_2}~\cblu{x\sub{0001}~x\sub{0100}~\beta\sub{0101}~x\sub{1000}~\beta\sub{1001}~\beta\sub{1100}})}\\
  &{(\cblu{\ff\sub{1110}} :
  \cred{A_2}~\cblu{x\sub{0010}~x\sub{0100}~\beta\sub{0110}~x\sub{1000}~\beta\sub{1010}~\beta\sub{1100}})} \to
  \Type
\end{align*}
First, we have a type $\cred{A_0}$ of \emph{points}. Second, we
have for every two points $\cblu{x\sub{01}},\cblu{x\sub{10}} : \cred{A_0}$, a type
$\cred{A_1}~\cblu{x\sub{01}~x\sub{10}}$ of \emph{lines} joining $\cblu{x\sub{01}}$ and $\cblu{x\sub{10}}$.
Third, for every three points $\cblu{x\sub{001}},\cblu{x\sub{010}},\cblu{x\sub{100}} : \cred{A_0}$ and
three lines $\cblu{\beta\sub{011}} : \cred{A_1}~\cblu{x\sub{001}~x\sub{010}}$, $\cblu{\beta\sub{101}} :
\cred{A_1}~\cblu{x\sub{001}~x\sub{100}}$, and $\cblu{\beta\sub{110}} : \cred{A_1}~\cblu{x\sub{010}~x\sub{100}}$, a type
$\cred{A_2}~\cblu{x\sub{001}~x\sub{010}~\beta\sub{011}~x\sub{100}~\beta\sub{101}~\beta\sub{110}}$ of \emph{triangles} with the given
boundary. The pattern continues with \emph{tetrahedra}, which we may also, more
technically, call $\cred{3}$-simplices. In general, $\cred{A}_\n$ describes the type of
$\n$-simplices indexed by their boundaries.

The second, known as the \emph{fibred} formulation, forgets the indexing and assembles each of the $\n$-simplex types into a single total space. We thus have types $\cred{X}_\n$ for $\n \geq \cred{0}$. The information about the faces of an $\n$-simplex is then encoded by a collection of \emph{face maps} $\cred{\partial}_\k : \cred{X}_\n \to \cred{X}_{\n-\cred{1}}$, for $\k \in \set{\cred{0}, \ldots, \n}$. Now for example, if $\cblu{\ff} : \cred{X_2}$ is a triangle, then both $\cred{\partial_0}~\cred{\partial_0}~\cblu{\ff}$ and $\cred{\partial_0}~\cred{\partial_1}~\cblu{\ff}$ should denote the same boundary point, so we must impose a coherence on the face maps. This says that we have homotopies:
\[\cred{\alpha}_{\k,\,\l} : \cred{\partial}_\k \circ \cred{\partial}_\l \simeq \cred{\partial}_{\l-\cred{1}} \circ \cred{\partial}_\k \quad \text{for } \k < \l\]
However, since we are working in a homotopical setting, the $\cred{\alpha}_{\k,\,\l}$ carry data, so it's not \emph{`just the case that two things are equal, with that being the end of the story'}. Rather, for $\k < \l < \m$, we can prove that $\cred{\partial}_\k \circ \cred{\partial}_\l \circ \cred{\partial}_\m \simeq \cred{\partial}_{\m-\cred{2}} \circ \cred{\partial}_{\l-\cred{1}} \circ \cred{\partial}_\k$ in two different ways, as illustrated by the following diagram:
\[
\begin{tikzpicture}
\path[draw, postaction={on each segment={mid arrow}}]
(1,1) -- (2,2) -- (3.4,2) -- (4.4,1)
(1,1) -- (2,0) -- (3.4,0) -- (4.4,1);
\filldraw (1,1) circle (1.5pt);
\filldraw (2,2) circle (1.5pt);
\filldraw (3.4,2) circle (1.5pt);
\filldraw (4.4,1) circle (1.5pt);
\filldraw (2,0) circle (1.5pt);
\filldraw (3.4,0) circle (1.5pt);
\draw (1,1) node[anchor=east]{${\cred{\partial}_\cblu{k} \circ \cred{\partial}_\cblu{l} \circ \cred{\partial}_\cblu{m}}$};
\draw (4.4,1) node[anchor=west]{${\cred{\partial}_{\m-\cred{2}} \circ \cred{\partial}_{\l-\cred{1}} \circ \cred{\partial}_\cblu{k}}$};
\draw (2,2) node[anchor=south east]{${\cred{\partial}_{\l-\cred{1}} \circ \cred{\partial}_\cblu{k} \circ \cred{\partial}_\cblu{m}}$};
\draw (3.4,2) node[anchor=south west]{${\cred{\partial}_{\l-\cred{1}} \circ \cred{\partial}_{\m-\cred{1}} \circ \cred{\partial}_\cblu{k}}$};
\draw (2,0) node[anchor=north east]{${\cred{\partial}_\cblu{k} \circ \cred{\partial}_{\m-\cred{1}} \circ \cred{\partial}_\cblu{l}}$};
\draw (3.4,0) node[anchor=north west]{${\cred{\partial}_{\m-\cred{2}} \circ \cred{\partial}_\cblu{k} \circ \cred{\partial}_\cblu{l}}$};
\draw (2.7,1) node{$\Downarrow$};
\end{tikzpicture}
\]
We may thus consider imposing higher coherences of the form:
\[\cred{\beta}_{\k,\,\l,\,\m} : \cred{\alpha}_{\k,\,\l} \star \cred{\partial}_\m \dot \cred{\partial}_{\l-\cred{1}} \star \cred{\alpha}_{\k,\,\m} \dot \cred{\alpha}_{\l-\cred{1},\,\m-\cred{1}} \star \cred{\partial}_\k
\simeq \cred{\partial}_\k \star \cred{\alpha}_{\l,\,\m} \dot \cred{\alpha}_{\k,\,\m-\cred{1}} \star \cred{\partial}_\l \dot \cred{\partial}_{\m-\cred{2}} \star \cred{\alpha}_{\k,\,\l}\]
As it turns out, in the truncated cases that consider only the first few simplex types, we can show that the indexed and fibred formulations are equivalent \emph{only when higher coherences are included in the latter}.~\cite{blog1} Thus, if we would like, in homotopy theory, the ability to talk about an $\n$-simplex with given boundary, the coherences are here to stay.

Now, of course the issue with these coherences persists; because the $\cred{\beta}_{\k,\,\l,\,\m}$ consist of data, we have to impose even higher coherences. At the next dimension, this is illustrated as follows, borrowing a diagram by Tilman Piesk~\cite{permute}:
\[
\includegraphics[scale=0.045]{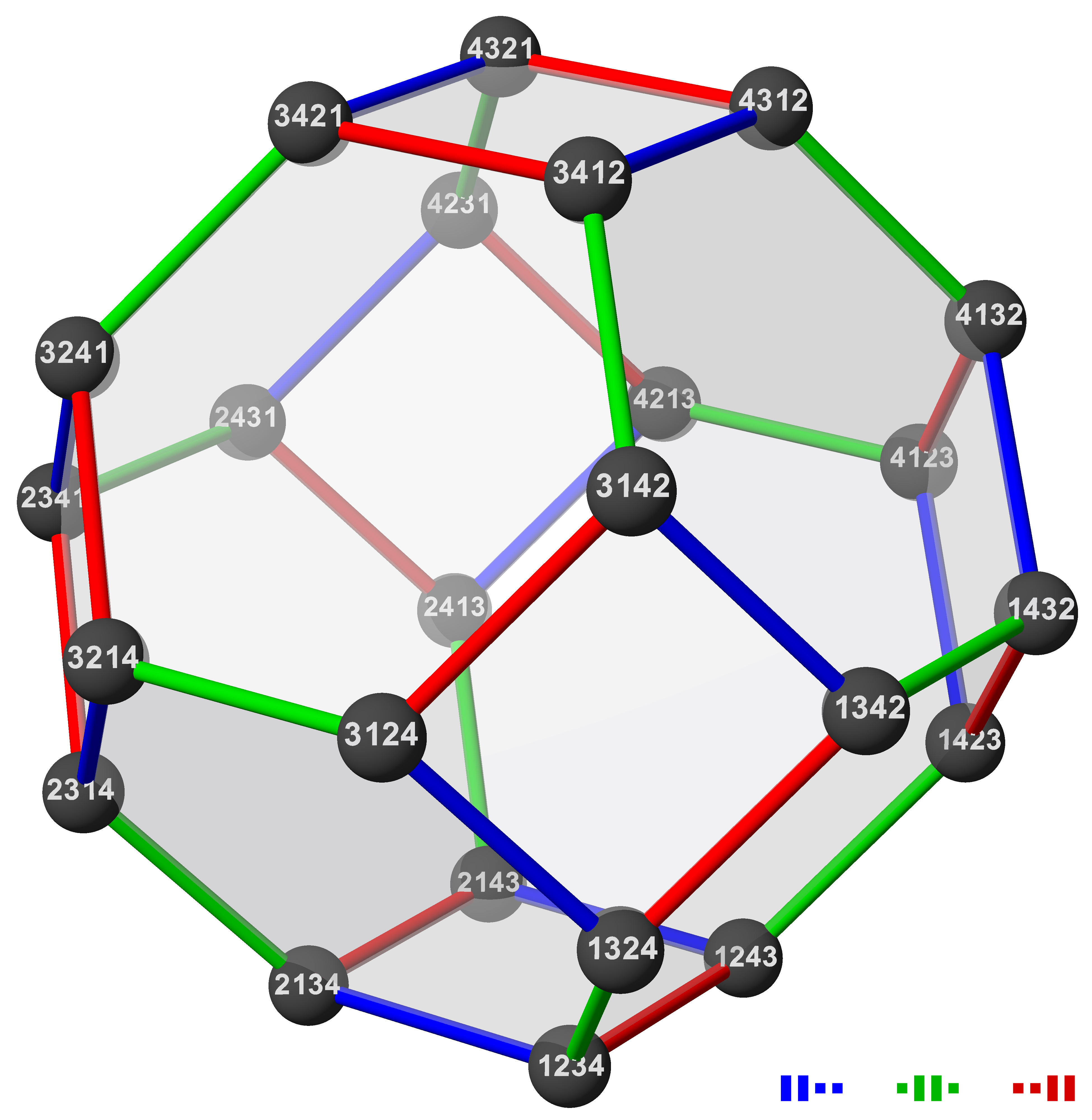} 
\]
In the context of coherence realtions, the hexagon is known as the $\cred{3}$-permutahedron and the truncated octahedron is known as the $\cred{4}$-permutahedron. In general, there are $\n$-permutahedra for every $\n \geq \cred{1}$. \emph{Our goal in this paper is to write down formulas for the coherence relations corresponding to them in any dimension.}

\noindent
\textbf{The Setting:} We will work in a different setting than the one described above. Let $\X$ be a topological space and let $\c : \X \to \X \to \X$ be a multiplication operation. We will assume that $\cblu{\mu}$ is \emph{strictly associative}, i.e. on the nose and not just up to homotopy. Our goal is to define what it means for $\c$ to be commutative up to homotopy in an infinitely coherent way. \bbox



\begin{center}
{\large \textbf{Part I: Combinatorics}}
\end{center}

There are $\n \cred{!}$ permutations on $\n$-elements, but what is further true is that there is a direct encoding of $\n$-permutations as terms of the type $\mathsf{\cred{Fin}}~\cred{1} \times \mathsf{\cred{Fin}}~\cred{2} \times \cdots \times \mathsf{\cred{Fin}}~\n$. We take the permutations to be on the first $\n$ capital letters of the Roman alphabet. An example of a $\cred{5}$-permutation is then given by $\sqbkt{\E}\,\sqbkt{\B}\,\sqbkt{\D}\,\sqbkt{\A}\,\sqbkt{\C}$. To obtain the desired encoding, we consider the permutations induced on the last $\k$ letters:
\begin{align*}
\dsqbkt{\E}& \mapsto \cred{0} \\
\sqbkt{\E}\,\dsqbkt{\D}& \mapsto \cred{1} \\
\sqbkt{\E}\,\sqbkt{\D}\,\dsqbkt{\C}& \mapsto \cred{2} \\
\sqbkt{\E}\,\dsqbkt{\B}\,\sqbkt{\D}\,\sqbkt{\C}& \mapsto \cred{1} \\
\sqbkt{\E}\,\sqbkt{\B}\,\sqbkt{\D}\,\dsqbkt{\A}\,\sqbkt{\C}& \mapsto \cred{3}
\end{align*}
When we add the $\k$-th last letter, we count its position in an zero-indexed manner from the left. The permutation $\sqbkt{\E}\,\sqbkt{\B}\,\sqbkt{\D}\,\sqbkt{\A}\,\sqbkt{\C}$ is thus sent to $(\cred{0}, \cred{1},\cred{2},\cred{1},\cred{3})$.

The facets of a $(\n+\cred{1})$-permutahedron are described by the ordered partitions of the first $\n+\cred{1}$ letters into non-empty groups that are internally unordered.~\cite{pons2023combinatorics} For example, two facets of the $\cred{5}$-permutahedron are given by $\sqbkt{\B\D}\,\sqbkt{\C}\,\sqbkt{\A\E}$ and $\sqbkt{\B\C}\,\sqbkt{\A\D\E}$; note that within each group the letters are alphabetically ordered. We can describe an algorithm, similar to the above, that extracts from each facet a region of the cube $\sqbkt{\cred{0},\cred{0}} \times \sqbkt{\cred{0},\cred{1}} \times \sqbkt{\cred{0},\cred{2}} \times \cdots \times \sqbkt{\cred{0},\n}$, such that all such regions disjointly partition the whole cube. In the cases of the above examples, we would get:
\begin{align*}
\dsqbkt{\E}& \mapsto \set{\cred{0}} &
\dsqbkt{\E}& \mapsto \set{\cred{0}}\\
\dsqbkt{\D}\,\sqbkt{\E}& \mapsto \set{\cred{0}} &
\dsqbkt{\D\E}& \mapsto \bkt{\cred{0},\cred{1}} \\
\sqbkt{\D}\,\dsqbkt{\C}\,\sqbkt{\E}& \mapsto \set{\cred{1}} &
\dsqbkt{\C}\,\sqbkt{\D\E}& \mapsto \set{\cred{0}}\\
\dsqbkt{\B\D}\,\sqbkt{\C}\,\sqbkt{\E}& \mapsto \bkt{\cred{0}, \cred{1}} &
\dsqbkt{\B\C}\,\sqbkt{\D\E}& \mapsto \bkt{\cred{0},\cred{1}} \\
\sqbkt{\B\D}\,\sqbkt{\C}\,\dsqbkt{\A\E}& \mapsto \bkt{\cred{3}, \cred{4}} &
\sqbkt{\B\C}\,\dsqbkt{\A\D\E}& \mapsto \bkt{\cred{2},\cred{4}} 
\end{align*}
We mark, with double square brackets, the \emph{`active'} group to which the latest letter was added. We think of a group as an uncertain ordering of its constituent letters. Thus, for example, given $\dsqbkt{\D\E}$, we extract not an index, but rather a range $\bkt{\cred{0},\cred{1}}$ for possible positions of the latest letter $\D$. In general, this range is not discrete, as we think of the act of permuting letters as a kind of continuous deformation. In order for the regions to be disjoint, we take the range to be open, which corresponds to thinking of $\dsqbkt{\D\E}$ as consisting of strictly intermediate stages of the deformation. With this in mind, the two examples above would be sent to the regions $\set{\cred{0}} \times \set{\cred{0}} \times \set{\cred{1}} \times \bkt{\cred{0},\cred{1}} \times \bkt{\cred{3},\cred{4}}$ and $\set{\cred{0}} \times \bkt{\cred{0},\cred{1}} \times \set{\cred{0}} \times \bkt{\cred{0},\cred{1}} \times \bkt{\cred{2},\cred{4}}$, respectively.

In what follows, we modify this algorithm slightly; since this process always begins with $\dsqbkt{\blank} \mapsto \set{\cred{0}}$, for the last letter $\bkt{\blank}$, we drop this singleton range and omit the $\bkt{\sqbkt{\cred{0},\cred{0}} \times \blank}$ from the cube. We thus consider the resulting regions to lie in the cube $\sqbkt{\cred{0},\cred{1}} \times \sqbkt{\cred{0},\cred{2}} \times \cdots \times \sqbkt{\cred{0},\n}$, which is the standard $\n$-cube of dimensions $\cred{1} \times \cred{2} \times \cdots \times \n$. Note that the top dimensional facet consisting of a single group is always sent to the whole interior of the cube.

\begin{figure}[h]
\vspace*{-1.5cm}
\noindent
\makebox[\textwidth]{
\scalebox{0.69}{
\begin{tikzpicture}
\draw (2.5,7.5,0) node[fill=white]{[\A\B]\,[\C\D]};
\draw (10,7.5,0) node[fill=white]{[\B]\,[\A\C\D]};
\draw (5,2.5,0) node[fill=white]{[\A\B\C]\,[\D]};
\draw (12.5,2.5,0) node[fill=white]{[\B\C]\,[\A\D]};
\draw (-2.5,0,0) node[fill=white]{[\A]\,[\B\C\D]};
\draw (17.5,0,0) node[fill=white]{[\B\C\D]\,[\A]};
\draw (2.5,-2.5,0) node[fill=white]{[\A\C]\,[\B\D]};
\draw (10,-2.5,0) node[fill=white]{[\C]\,[\A\B\D]};
\draw (5,-7.5,0) node[fill=white]{[\A\C\D]\,[\B]};
\draw (12.5,-7.5,0) node[fill=white]{[\C\D]\,[\A\B]};
\draw (10,-12.5,0) node[fill=white]{[\D]\,[\A\B\C]};
\draw (2.5,-12.5,0) node[fill=white]{[\A\D]\,[\B\C]};
\draw (5,-17.5,0) node[fill=white]{[\A\B\D]\,[\C]};
\draw (12.5,-17.5,0) node[fill=white]{[\B\D]\,[\A\C]};
\path[draw]
(5,10,0) -- (10,10,0) node[midway,anchor=north]{[\B]\,[\A\D]\,[\C]}
(10,10,0) -- (15,10,0) node[midway,anchor=north]{[\B]\,[\D]\,[\A\C]}
(5,10,0) -- (5,5,0) node[midway,anchor=west]{[\B]\,[\A]\,[\C\D]}
(15,5,0) -- (15,10,0) node[midway,anchor=east]{[\B]\,[\C\D]\,[\A]}
(-5,5,0) -- (0,5,0) node[midway,anchor=north]{[\A]\,[\B]\,[\C\D]}
(0,5,0) -- (5,5,0) node[midway,anchor=north]{[\A\B]\,[\C]\,[\D]}
(5,5,0) -- (10,5,0) node[midway,anchor=south]{[\B]\,[\A\C]\,[\D]}
(10,5,0) -- (15,5,0) node[midway,anchor=south]{[\B]\,[\C]\,[\A\D]}
(-5,0,0) -- (-5,5,0) node[midway,anchor=west]{[\A]\,[\B\D]\,[\C]}
(0,0,0) -- (0,5,0) node[midway,anchor=east]{[\A]\,[\B\C]\,[\D]}
(10,0,0) -- (10,5,0) node[midway,anchor=east]{[\B\C]\,[\A]\,[\D]}
(15,0,0) -- (15,5,0) node[midway,anchor=west]{[\B\C]\,[\D]\,[\A]}
(20,0,0) -- (20,5,0) node[midway,anchor=east]{[\B\D]\,[\C]\,[\A]}
(0,0,0) -- (5,0,0) node[midway,anchor=south]{[\A\C]\,[\B]\,[\D]}
(5,0,0) -- (10,0,0) node[midway,anchor=north]{[\C]\,[\A\B]\,[\D]}
(10,0,0) -- (15,0,0) node[midway,anchor=north]{[\C]\,[\B]\,[\A\D]}
(-5,0,0) -- (-5,-5,0) node[midway,anchor=west]{[\A]\,[\D]\,[\B\C]}
(0,0,0) -- (0,-5,0) node[midway,anchor=east]{[\A]\,[\C]\,[\B\D]}
(5,0,0) -- (5,-5,0) node[midway,anchor=west]{[\C]\,[\A]\,[\B\D]}
(15,0,0) -- (15,-5,0) node[midway,anchor=east]{[\C]\,[\B\D]\,[\A]}
(-5,-5,0) -- (0,-5,0) node[midway,anchor=south]{[\A]\,[\C\D]\,[\B]}
(0,-5,0) -- (5,-5,0) node[midway,anchor=north]{[\A\C]\,[\D]\,[\B]}
(5,-5,0) -- (10,-5,0) node[midway,anchor=south]{[\C]\,[\A\D]\,[\B]}
(10,-5,0) -- (15,-5,0) node[midway,anchor=south]{[\C]\,[\D]\,[\A\B]}
(15,-5,0) -- (20,-5,0) node[midway,anchor=south]{[\C\D]\,[\B]\,[\A]}
(10,-5,0) -- (10,-10,0) node[midway,anchor=east]{[\C\D]\,[\A]\,[\B]}
(0,-10,0) -- (5,-10,0) node[midway,anchor=south]{[\A\D]\,[\C]\,[\B]}
(5,-10,0) -- (10,-10,0) node[midway,anchor=north]{[\D]\,[\A\C]\,[\B]}
(10,-10,0) -- (15,-10,0) node[midway,anchor=north]{[\D]\,[\C]\,[\A\B]}
(5,-10,0) -- (5,-15,0) node[midway,anchor=west]{[\D]\,[\A]\,[\B\C]}
(15,-10,0) -- (15,-15,0) node[midway,anchor=east]{[\D]\,[\B\C]\,[\A]}
(0,-15,0) -- (5,-15,0) node[midway,anchor=north]{[\A\D]\,[\B]\,[\C]}
(5,-15,0) -- (10,-15,0) node[midway,anchor=south]{[\D]\,[\A\B]\,[\C]}
(10,-15,0) -- (15,-15,0) node[midway,anchor=south]{[\D]\,[\B]\,[\A\C]}
(10,-15,0) -- (10,-20,0) node[midway,anchor=east]{[\B\D]\,[\A]\,[\C]}
(0,-20,0) -- (5,-20,0) node[midway,anchor=south]{[\A\B]\,[\D]\,[\C]};
\path[draw, ultra thick, color=darkgreen]
(0,-5,0) -- (0,5,0) -- (15,5,0) -- (15,0,0) -- (15,-5,0) -- (0,-5,0)
(0,-5,0)--(-5,-5,0)--(-5,5,0)--(0,5,0)
(15,-5,0)--(20,-5,0)--(20,5,0)--(15,5,0)
(0,5,0) -- (0,10,0) -- (15,10,0)--(15,5,0)
(0,-5,0) -- (0,-10,0) -- (15,-10,0)--(15,-5,0)
(0,-10,0) -- (0,-20,0) -- (15,-20,0)--(15,-10,0);
\draw (0,9.84,0) node[fill=white]{[\A]\,[\B]\,[\D]\,[\C]};
\draw (5,9.84,0) node[fill=white]{[\B]\,[\A]\,[\D]\,[\C]};
\draw (0,4.84,0) node[fill=white]{[\A]\,[\B]\,[\C]\,[\D]};
\draw (5,5,0) node[fill=white]{[\B]\,[\A]\,[\C]\,[\D]};
\draw (10,5,0) node[fill=white]{[\B]\,[\C]\,[\A]\,[\D]};
\draw (15,4.84,0) node[fill=white]{[\B]\,[\C]\,[\D]\,[\A]};
\draw (0,0,0) node[fill=white]{[\A]\,[\C]\,[\B]\,[\D]};
\draw (5,0,0) node[fill=white]{[\C]\,[\A]\,[\B]\,[\D]};
\draw (10,0,0) node[fill=white]{[\C]\,[\B]\,[\A]\,[\D]};
\draw (15,0,0) node[fill=white]{[\C]\,[\B]\,[\D]\,[\A]};
\draw (0,-4.84,0) node[fill=white]{[\A]\,[\C]\,[\D]\,[\B]};
\draw (5,-5,0) node[fill=white]{[\C]\,[\A]\,[\D]\,[\B]};
\draw (10,-5,0) node[fill=white]{[\C]\,[\D]\,[\A]\,[\B]};
\draw (15,-4.84,0) node[fill=white]{[\C]\,[\D]\,[\B]\,[\A]};
\draw (0,-10.16,0) node[fill=white]{[\A]\,[\D]\,[\C]\,[\B]};
\draw (5,-10,0) node[fill=white]{[\D]\,[\A]\,[\C]\,[\B]};
\draw (10,-10,0) node[fill=white]{[\D]\,[\C]\,[\A]\,[\B]};
\draw (15,-10.16,0) node[fill=white]{[\D]\,[\C]\,[\B]\,[\A]};
\draw (0,-15,0) node[fill=white]{[\A]\,[\D]\,[\B]\,[\C]};
\draw (5,-15,0) node[fill=white]{[\D]\,[\A]\,[\B]\,[\C]};
\draw (10,-15,0) node[fill=white]{[\D]\,[\B]\,[\A]\,[\C]};
\draw (15,-15,0) node[fill=white]{[\D]\,[\B]\,[\C]\,[\A]};
\draw (10,-19.84,0) node[fill=white]{[\B]\,[\D]\,[\A]\,[\C]};
\draw (15,-19.84,0) node[fill=white]{[\B]\,[\D]\,[\C]\,[\A]};
\end{tikzpicture}}}
\caption{The net of the $\cred{4}$-permutahedron}
\label{fig:net}
\end{figure}
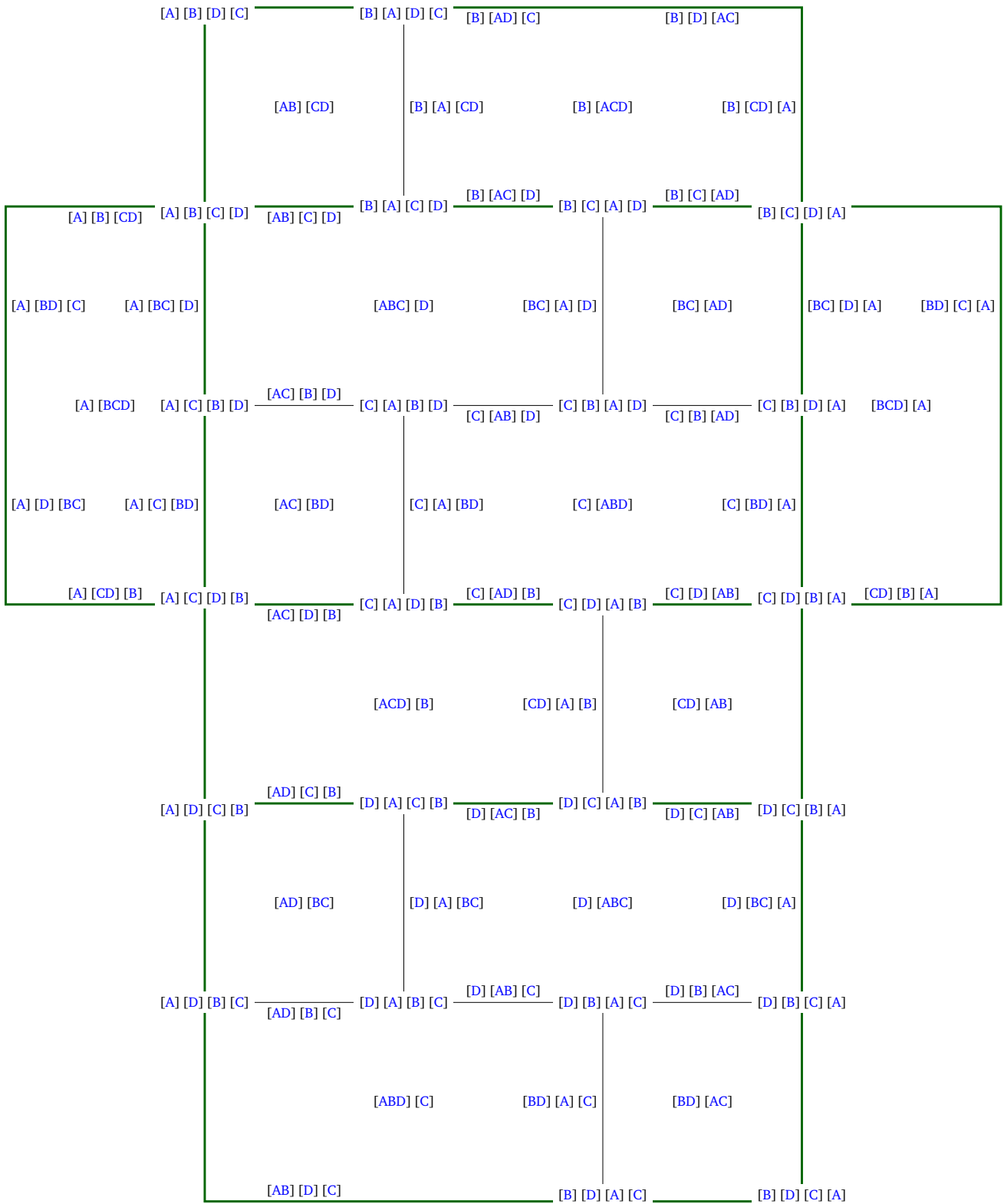

According to the above placement algorithm the facets of the $\cred{3}$-permutahedron are laid out as follows:
\[
\begin{tikzpicture}
\path[draw]
(0,0) -- (0,2) node[midway,anchor=east]{[\A]\,[\B\C]}
(0,2) -- (4,2) node[midway,anchor=south]{[\A\B]\,[\C]}
(4,2) -- (8,2) node[midway,anchor=south]{[\B]\,[\A\C]}
(0,0) -- (4,0) node[midway,anchor=north]{[\A\C]\,[\B]}
(4,0) -- (8,0) node[midway,anchor=north]{[\C]\,[\A\B]}
(8,0) -- (8,2) node[midway,anchor=west]{[\B\C]\,[\A]};
\draw (0,2) node[fill=white]{[\A]\,[\B]\,[\C]};
\draw (0,0) node[fill=white]{[\A]\,[\C]\,[\B]};
\draw (4,2) node[fill=white]{[\B]\,[\A]\,[\C]};
\draw (8,2) node[fill=white]{[\B]\,[\C]\,[\A]};
\draw (4,0) node[fill=white]{[\C]\,[\A]\,[\B]};
\draw (8,0) node[fill=white]{[\C]\,[\B]\,[\A]};
\draw (4,1) node[fill=white]{[\A\B\C]};
\end{tikzpicture}
\]
In the case of the $\cred{4}$-permutahedron, we represent the placements in a net model given in Figure \ref{fig:net} on the following page. This model may be physically cut out and assembled into a $\cred{1} \times \cred{2} \times \cred{3}$ cube, using the overhanging tabs to glue.
\vspace{0.2cm}

\noindent
\textbf{Theorem}: Every point $\p : \sqbkt{\cred{0},\cred{1}} \times \sqbkt{\cred{0},\cred{2}} \times \cdots \times \sqbkt{\cred{0},\n}$ lies in a unique facet of the $(\n+\cred{1})$-permutahedron, and thus the facets disjointly partition the cube.
\vspace{0.2cm}

\noindent 
\textbf{Idea}: Given $(\cred{1}, \cred{1.7}, \cred{2}, \cred{0.5}) :  \sqbkt{\cred{0},\cred{1}} \times \sqbkt{\cred{0},\cred{2}} \times \sqbkt{\cred{0},\cred{3}} \times \sqbkt{\cred{0},\cred{4}}$, we find the enclosing facet as follows:
\[
\begin{tikzpicture}
\draw (0,0) node[anchor=west]{[\E]};
\draw[(-), thick, color=blue] (-0.1,-0.35) -- (-0.75,-0.35);
\draw (-0.1,-0.35) node[anchor=south]{\cred{1}};
\draw (-0.75,-0.35) node[anchor=south]{\cred{0}};
\filldraw[color=red] (-0.1,-0.35) circle (1.5pt);
\draw (0,-0.7) node[anchor=west]{[\E]\,[\D]};
\draw[(-), thick, color=blue]
(-0.1,-1.05) -- (-0.75,-1.05);
\draw[(-), thick, color=blue]
(-0.75,-1.05)--(-1.4,-1.05);
\draw (-0.1,-1.05) node[anchor=south]{\cred{2}};
\draw (-0.75,-1.05) node[anchor=south]{\cred{1}};
\draw (-1.4,-1.05) node[anchor=south]{\cred{0}};
\filldraw[color=red] (-0.295,-1.05) circle (1.5pt);
\draw (0,-1.4) node[anchor=west]{[\E]\,[\C\D]};
\draw[(-), thick, color=blue]
(-0.1,-1.75) -- (-1.4,-1.75);
\draw[(-), thick, color=blue]
(-1.4,-1.75) -- (-2.05,-1.75);
\draw (-0.1,-1.75) node[anchor=south]{\cred{3}};
\draw (-1.4,-1.75) node[anchor=south]{\cred{1}};
\draw (-2.05,-1.75) node[anchor=south]{\cred{0}};
\filldraw[color=red] (-0.75,-1.75) circle (1.5pt);
\draw (0,-2.1) node[anchor=west]{[\E]\,[\B\C\D]};
\draw[(-), thick, color=blue]
(-0.1,-2.45) -- (-2.05,-2.45);
\draw[(-), thick, color=blue]
(-2.05,-2.45) -- (-2.70,-2.45);
\draw (-0.1,-2.45) node[anchor=south]{\cred{4}};
\draw (-2.05,-2.45) node[anchor=south]{\cred{1}};
\draw (-2.70,-2.45) node[anchor=south]{\cred{0}};
\filldraw[color=red] (-2.375,-2.45) circle (1.5pt);
\draw (0,-2.8) node[anchor=west]{[\A\E]\,[\B\C\D]};
\end{tikzpicture}
\]
\vspace{-0.7cm}

\noindent
\textbf{Proof}: For a facet $\F$, let $\F_\k$ denote its restriction to the last $\k$ letters. According to our algorithm for situating a facet in space, the first $\k$ coordinates of the region corresponding to $\F$ are determined by $\F_{\k+\cred{1}}$. Now let $\p : \sqbkt{\cred{0},\cred{1}} \times \sqbkt{\cred{0},\cred{2}} \times \cdots \times \sqbkt{\cred{0},\n}$ be our point. Let $\p_\k$ denote its projection to the first $\k$ coordinates for $\k \geq \cred{0}$, with $\p_\cred{0} = \bkt{}$. The point $\p_\cred{0}$ is said to lie in the unique facet $[\blank]$ of the $\cred{1}$-permutahedron. Now suppose that $\n = \k + \cred{1}$. By induction, $\p_\k$ lies in the region of a unique facet $\F$ of the $\bkt{\k+\cred{1}}$-pemutahedron. Note that, if $\p$ lies in the facet $\G$ of the $\bkt{\k+\cred{2}}$-pemutahedron, then $\G_{\k+\cred{1}} = \F$, so we only have to consider facets $\G$ that extend $\F$ by the first letter.

Reading from left to right, suppose that $\F$ consists of $\m$ groups of size $\bkt{\cblu{\alpha}_\cred{1}, \cdots, \cblu{\alpha}_\m}$. We have that $\sum_{\j = \cred{1}}^{\m} \cblu{\alpha}_\j= \k + \cred{1} = \n$, and the possible positions that the new letter can take also continuously range in $\sqbkt{\cred{0},\n}$. If we add the new letter as a discrete group directly following group $\l$ for $\l \geq \cred{0}$, then the situating algorithm takes $\G$ to $\set{\sum_{\j = \cred{1}}^{\l} \cblu{\alpha}_\j}$ in the last coordinate. On the other hand, if we add the new letter to the group $\l$, then then the situating algorithm takes $\G$ to $\bkt{\sum_{\j = \cred{1}}^{\l-\cred{1}} \cblu{\alpha}_\j, \sum_{\j = \cred{1}}^{\l} \cblu{\alpha}_\j}$ in the last coordinate. Note that these regions disjointly partition $\sqbkt{\cred{0},\n}$, thus there is a unique $\G$ such that $\p$ lies in its region. \qed
\vspace{0.2cm}

Note that in the diagram for the \emph{Idea} above, we draw the partition of $\sqbkt{\cred{0},\n}$ that arises in the proof of the theorem, and mark where the latest coordinate lies. This tells us where to add the latest letter. \bbox

\clearpage

\clearpage

\begin{center}
{\cblu{\circled{\large \cred{\quarterrest}}}\ \ {\large \textbf{Interlude}}}
\end{center}

From the above, we know how to rigidly transform the $(\n+\cred{1})$-permutahedron into an $\n$-cube. In this interlude, we discuss the forms that our equations corresponding to these cubes will take.

Note that a decomposition of the boundary of a cube does not automatically give rise to conveniently expressible equations. For example, in Figure \ref{fig:net}, we saw that the $\cred{2} \times \cred{3}$ faces looked like the following diagram on the left, and we can write a $\cred{2}$-categorical composition of these rectangles by first composing each column vertically, and then composing the results horizontally. On the other hand, such a convenient composition is not possible to express in the following \emph{tatami mat} diagram to the right.
\[
\begin{tikzpicture}
\path[draw, color=darkgreen]
(0,0) -- (0,2.25) -- (1.5,2.25) -- (1.5,0) -- (0,0)
(0.75,0) -- (0.75,2.25)
(0,0.75)--(0.75,0.75)
(0.75,1.5)--(1.5,1.5);
\end{tikzpicture}
\quad \quad \quad \quad
\begin{tikzpicture}
\path[draw, color=darkgreen]
(0,0) -- (2.25,0) -- (2.25,2.25) -- (0,2.25) -- (0,0)
(0.75,0)--(0.75,1.5) -- (0,1.5)
(0,1.5) -- (1.5,1.5) -- (1.5,2.25)
(1.5,2.25) -- (1.5,0.75) -- (2.25,0.75)
(2.25,0.75) -- (0.75,0.75);
\end{tikzpicture}
\]
We shall see that the former case is representative, and that the general coherence equations can be expressed as cube filling problems in which the faces are expressed as iterated compositions along successive axes.

Note that along the middle vertical axis in the rectangle, the composite of the lines on the left would be something like $\bkt{[\A\C]\,[\B]\,[\D] \dot [\C]\,[\A\B]\,[\D]} \dot [\C]\,[\B]\,[\A\D]$ whereas on the right, it would be something like $[\A\C]\,[\D]\,[\B] \dot \bkt{[\C]\,[\A\B]\,[\D] \dot [\C]\,[\B]\,[\A\D]}$. Unfortunately, in many settings, path composition is only be weakly associative, and putative formulas would have to account for reassociation. \emph{We would like to avoid this.} For this reason, we work in a particularly nice setting for dealing with (higher) homotopies. This will be the setting of \emph{cubical Moore paths}.

Let $\X$ be a space, $\a,\b : \X$, and $\n : \cred{\mathbb{N}}$, then we define $\Path_\n~\X~\a~\b$ to be the subspace of $\X^{\sqbkt{\cred{0},\n}}$ consisting of such $\cblu{\alpha} : \X^{\sqbkt{\cred{0},\n}}$ with $\cblu{\alpha}(\cred{0}) \equiv \a$ and $\cblu{\alpha}(\n) \equiv \b$. There is an evident composition operation $\dot : \Path_\n~\X~\a~\b \to \Path_\m~\X~\b~\rc \to \Path_{\n+\m}~\X~\a~\rc$ obtained by concatenation; note that this operation is definitionally associative.

Now suppose that $\Y : \sqbkt{\cred{0},\n} \to \Top$ is a family of topological spaces parameterised over the interval $\sqbkt{\cred{0},\n}$ continuously, and that $\a : \Y(\cred{0})$ and $\b : \Y(\n)$, then define $\PathP_\n~\Y~\a~\b$ to be the subspace of dependent functions $\cblu{\alpha} : (\i : \sqbkt{\cred{0},\n}) \to \Y~\i$ with $\cblu{\alpha}(\cred{0}) \equiv \a$ and $\cblu{\alpha}(\n) \equiv \b$. If $\Y_\cblu{0}$ and $\Y_\cblu{1}$ are composable lines of spaces, then there is a dependent composition operation $\dot : \PathP_\n~\Y_\cblu{0}~\a~\b \to \Path_\m~\Y_\cblu{1}~\b~\rc \to \PathP_{\n+\m}~(\Y_\cblu{0} \dot \Y_\cblu{1})~\a~\rc$, which is likewise definitionally associative.

If $\c : \X \to \X \to \X$ is a composition operation, on $\X$, then for points $\a, \b, \cblu{c} : \X$ and a path $\cblu{\alpha} : \Path_\n~\X~\b~\cblu{c} $, we have that $\cred{\lambda}~\i\text{.}~ \c~\a~\bkt{\cblu{\alpha}~\i} : \Path_\n~\X~(\c~\a~\b)~(\c~\a~\cblu{c})$. This is a whiskering operation that we can define in terms of our primitive operations. Note that whiskering distributes over composition definitionally, such that:
\[\cred{\lambda}~\i\text{.}~ \c~\a~\bkt{\bkt{\cblu{\alpha} \dot \cblu{\beta}}~\i} \equiv \bkt{\cred{\lambda}~\i\text{.}~ \c~\a~\bkt{\cblu{\alpha}~\i}} \dot \bkt{{\cred{\lambda}~\i\text{.}~ \c~\a~\bkt{\cblu{\beta}~\i}}}\]
\clearpage

Using these operations, we are able to define the first few levels of the permutahedral coherence relations. We have:
\begin{align*}
\bkt{\Perm~\X~\c}_\cred{0} &\equiv \X \\
\bkt{\Perm~\X~\c}_\cred{1}~\A~\B &\equiv \Path_\cred{1}~\X~\bkt{\c~\A~\B}~\bkt{\c~\B~\A}
\end{align*}
At the next level, we iterate path types and form a $\cred{1} \times \cred{2}$ rectangle as follows:
\begin{align*}
&\bkt{\Perm~\X~\c}_\cred{2}~\A~\B~\A\B~\C~\A\C~\B\C \equiv \\
&\quad \PathP_\cred{1}~\bkt{\cred{\lambda}~\i\text{.}~ \Path_\cred{2}~\X~\bkt{\c~\A~\bkt{\B\C~\i}}~\bkt{\c~\bkt{\B\C~\i}~\A}} \\
&\quad\quad \bkt{\bkt{\cred{\lambda}~\j\text{.}~ \c~\bkt{\A\B~\j}~\C} \dot \bkt{\cred{\lambda}~\j\text{.}~ \c~\B~\bkt{\A\C~\j}}} \\
&\quad\quad \bkt{\bkt{\cred{\lambda}~\j\text{.}~ \c~\bkt{\A\C~\j}~\B} \dot \bkt{\cred{\lambda}~\j\text{.}~ \c~\C~\bkt{\A\B~\j}}}
\end{align*}
Note that this is well-typed only because $\c$ is definitionally associative. In the next section we will define the general $\bkt{\n+\cred{1}}$-st permutahedral coherence relation to be of the form:
\[\PathP_\cred{1}~\bkt{\cred{\lambda}~\i_\cblu{1}\text{.}~ \PathP_\cred{2}~\bkt{\cred{\lambda}~\i_\cblu{2}\text{.}~ \cdots~\bkt{\PathP_{\cblu{\n}-\cred{1}}~\bkt{\cred{\lambda}~\i_{\cblu{\n-1}}\text{.}~ \Path_\n~\X~\cprime{?} ~\cprime{?}}~\cprime{?}~\cprime{?} }~\cdots~}~\cprime{?} ~\cprime{?}}~\cprime{?} ~\cprime{?} \tag*{\bbox}\]

\begin{center}
{\large \textbf{Part II: Type Theory}}
\end{center}

\emph{Under the fictitious assumption} that every composition operation $\c : \X \to \X \to \X$ is definitionally associative, the $\SST$ of permutahedra is defined in the language of \emph{Displayed Type Theory (dTT)}~\cite{dTT} as follows:
\begin{lstlisting}
(*$\Perm$*) : ((*$\X$*) : (*$\Type$*)) ((*$\c$*) : (*$\X$*) (*$\to$*) (*$\X$*) (*$\to$*) (*$\X$*)) (*$\to$*) (*$\SST$*)
(*$\Z$*) ((*$\Perm$*) (*$\X$*) (*$\c$*)) = (*$\X$*)
(*$\S$*) ((*$\Perm$*) (*$\X$*) (*$\c$*)) (*$\A$*) =
  (*$\Perm\d$*)
    (*$\X$*) ((*$\lambda$*) (*$\B$*). (*$\Path$*) (*$\X$*) ((*$\c$*) (*$\A$*) (*$\B$*)) ((*$\c$*) (*$\B$*) (*$\A$*)))
    (*$\c$*) ((*$\lambda$*) (*$\B$*) (*$\A\B$*) (*$\C$*) (*$\A\C$*). ((*$\lambda$*) (*$\i$*). (*$\c$*) ((*$\A\B$*) (*$\i$*)) (*$\C$*)) (*$\dot$*) ((*$\lambda$*) (*$\i$*). (*$\c$*) (*$\B$*) ((*$\A\C$*) (*$\i$*))))
\end{lstlisting}
This definition is \emph{only well-typed up to associativity}; note that:
\begin{align*}
\Perm\d : &\bkt{\X : \Type}\;\bkt{\X\cblu{'} : \X \to \Type}\; \bkt{\c : \X \to \X \to \X} \\
&\bkt{\c\cblu{'} : \bkt{\x : \X} \to \X\cblu{'}\;\x \to \bkt{\y : \X} \to \X\cblu{'}\;\y \to \X\cblu{'}\;\bkt{\c\;\x\;\y}} \to \SST\d\;\bkt{\Perm\;\X\;\c}
\end{align*}
We define the case $\S\;\bkt{\Perm\;\X\;\c}\;\A$ by setting:
\begin{align*}
\X\cblu{'}\;\B &\equiv \Path\;\X\;\bkt{\c\;\A\;\B}\;\bkt{\c\;\B\;\A} \\
\c\cblu{'}\;\B\;\A\B\;\C\;\A\C &\equiv \bkt{\lambda\;\i\text{.}~ \c\;\bkt{\A\B\;\i}\;\C} \dot \bkt{\lambda\;\i\text{.}~\c\;\B\;\bkt{\A\C\;\i}}
\end{align*}
We have that:
\begin{align*}
\bkt{\lambda\;\i\text{.}~ \c\;\bkt{\A\B\;\i}\;\C} &: \Path\;\X\;\bkt{\c\;\bkt{\c\;\A\;\B}\;\C}\;\bkt{\c\;\bkt{\c\;\B\;\A}\;\C} \\
\bkt{\lambda\;\i\text{.}~\c\;\B\;\bkt{\A\C\;\i}} &: \Path\;\X\;\bkt{\c\;\B\;\bkt{\c\;\A\;\C}}\;\bkt{\c\;\B\;\bkt{\c\;\C\;\A}}
\end{align*}
Further, the expected type of $\c\cblu{'}\;\B\;\A\B\;\C\;\A\C$ is:
\[\X\cblu{'}\;\bkt{\c\;\B\;\C} \equiv \Path\;\X\;\bkt{\c\;\A\;\bkt{\c\;\B\;\C}}\;\bkt{\c\;\bkt{\c\;\B\;\C}\;\A}\]
So this definition of $\c\cblu{'}$ only makes sense if we definitionally have that:
\begin{align*}
\c\;\A\;\bkt{\c\;\B\;\C} &\equiv \c\;\bkt{\c\;\A\;\B}\;\C \\
\c\;\bkt{\c\;\B\;\A}\;\C &\equiv \c\;\B\;\bkt{\c\;\A\;\C} \\
\c\;\B\;\bkt{\c\;\C\;\A} &\equiv \c\;\bkt{\c\;\B\;\C}\;\A.
\end{align*}

One attempted remedy to this problem would be to add an associativity witness to the input data of $\Perm$ and thus to alternatively consider:
\begin{align*}
\Perm : &\bkt{\X : \Type}\;\bkt{\c : \X \to \X \to \X} \\
&\bkt{\a : \bkt{\x\;\y\;\z : \X} \to \Path\;\X\;\bkt{\c\;\x\;\bkt{\c\;\y\;\z}}\;\bkt{\c\;\bkt{\c\;\x\;\y}\;\z}} \to \SST  
\end{align*}
We would then be able to amend the definition of $\c\cblu{'}$ to read:
\begin{align*}
&\c\cblu{'}\;\B\;\A\B\;\C\;\A\C \\
&\quad \equiv \bkt{\a\;\A\;\B\;\C} \dot \bkt{\lambda\;\i\text{.}~ \c\;\bkt{\A\B\;\i}\;\C} \dot \bkt{\a\;\B\;\A\;\C}^\cred{-1} \dot \bkt{\lambda\;\i\text{.}~\c\;\B\;\bkt{\A\C\;\i}} \dot \bkt{\a\;\B\;\C\;\A}
\end{align*}
We are thus able to define a well-typed \emph{dependent composition operation}. The crux is then defining the \emph{dependent associator}, and here one runs into higher coherence issues. This is difficult to give a concise explanation of, but, in brief, we would want:
\begin{equation}
\begin{aligned}
&\a\cblu{'}\;\B\;\A\B\;\C\;\A\C\;\D\;\A\D \\
&\quad : \PathP\;\bkt{\lambda\;\i\text{.}~ \X\cblu{'}\;\bkt{\a\;\A\;\B\;\C\;\i}} \\
&\ \ \quad\quad \bkt{\c\cblu{'}\;\B\;\A\B\;\bkt{\c\;\C\;\D}\;\bkt{\c\cblu{'}\;\C\;\A\C\;\D\;\A\D}} \\
&\ \ \quad\quad \bkt{\c\cblu{'}\;\bkt{\c\;\B\;\C}\;\bkt{\c\cblu{'}\;\B\;\A\B\;\C\;\A\C}\; \D\;\A\D}
\end{aligned}
\tag{$\dagger$}
\label{eq:dep-assoc}
\end{equation}
We then have that each $\c\cblu{'}$ expands to the definition above, with three reassociations. The resulting expression is a $\PathP$ over associativity, with each endpoint having six associativity operations, resulting in a rectangle with fourteen association operations along its boundary. Writing this out carefully reveals that this can be decomposed to require the \emph{associahedral} coherence relations.

In our topological setting, the interpretation of this is that if $\c$ is not definitionally associative, but rather is only associative to homotopy, then we would require the data of that associativity proof to be infinitely coherent. \emph{This makes sense!} The main result of this paper is thus to solve the problem of defining the pattern of permutahedral coherence relations, \emph{up to the problem of defining the pattern of associahedral coherence relations}.

Let's investigate what we have. The first few $\n$-simplex types of $\Perm$ are as follows:
\begin{lstlisting}
((*$\Perm$*) (*$\X$*) (*$\c$*))(*$_\cred{0}$*) :=
  (*$\Z$*) ((*$\Perm$*) (*$\X$*) (*$\c$*)) =
  (*$\X$*)
((*$\Perm$*) (*$\X$*) (*$\c$*))(*$_\cred{1}$*) (*$\A$*) (*$\B$*) :=
  (*$\Z\d$*) ((*$\S$*) ((*$\Perm$*) (*$\X$*) (*$\c$*)) (*$\A$*)) (*$\B$*) =
  (*$\Path$*) (*$\X$*) ((*$\c$*) (*$\A$*) (*$\B$*)) ((*$\c$*) (*$\B$*) (*$\A$*))
((*$\Perm$*) (*$\X$*) (*$\c$*))(*$_\cred{2}$*) (*$\A$*) (*$\B$*) (*$\A\B$*) (*$\C$*) (*$\A\C$*) (*$\B\C$*) :=
  (*$\Z\d\d$*) ((*$\S\d$*) ((*$\S$*) ((*$\Perm$*) (*$\X$*) (*$\c$*)) (*$\A$*)) (*$\B$*) (*$\A\B$*)) (*$\C$*) (*$\A\C$*) (*$\B\C$*) =
  (*$\PathP$*)
    ((*$\lambda$*) (*$\i$*).
      (*$\Path$*) (*$\X$*) ((*$\c$*) (*$\A$*) ((*$\B\C$*) (*$\i$*))) ((*$\c$*) ((*$\B\C$*) (*$\i$*)) (*$\A$*)))
    (((*$\lambda$*) (*$\j$*). (*$\c$*) ((*$\A\B$*) (*$\j$*)) (*$\C$*)) (*$\dot$*) ((*$\lambda$*) (*$\j$*). (*$\c$*) (*$\B$*) ((*$\A\C$*) (*$\j$*))))
    (((*$\lambda$*) (*$\j$*). (*$\c$*) ((*$\A\C$*) (*$\j$*)) (*$\B$*)) (*$\dot$*) ((*$\lambda$*) (*$\j$*). (*$\c$*) (*$\C$*) ((*$\A\B$*) (*$\j$*))))
\end{lstlisting}
We now arrive at the coherence for the $\cred{4}$-permutahedron; one should verify that this is an equational representation of Figure \ref{fig:net}:
\begin{lstlisting}
((*$\Perm$*) (*$\X$*) (*$\c$*))(*$_\cred{3}$*) (*$\A$*) (*$\B$*) (*$\A\B$*) (*$\C$*) (*$\A\C$*) (*$\B\C$*) (*$\A\B\C$*) (*$\D$*) (*$\A\D$*) (*$\B\D$*) (*$\A\B\D$*) (*$\C\D$*) (*$\A\C\D$*) (*$\B\C\D$*)  :=
  (*$\Z\d\d\d$*) ((*$\S\d\d$*) ((*$\S\d$*) ((*$\S$*) ((*$\Perm$*) (*$\X$*) (*$\c$*)) (*$\A$*)) (*$\B$*) (*$\A\B$*)) (*$\C$*) (*$\A\C$*) (*$\B\C$*) (*$\A\B\C$*))
    (*$\D$*) (*$\A\D$*) (*$\B\D$*) (*$\A\B\D$*) (*$\C\D$*) (*$\A\C\D$*) (*$\B\C\D$*) =
  (*$\PathP$*)
    ((*$\lambda$*) (*$\i$*).
      (*$\PathP$*)
        ((*$\lambda$*) (*$\j$*).
          (*$\Path$*) (*$\X$*) ((*$\c$*) (*$\A$*) ((*$\B\C\D$*) (*$\i$*) (*$\j$*))) ((*$\c$*) ((*$\B\C\D$*) (*$\i$*) (*$\j$*)) (*$\A$*)))
        (((*$\lambda$*) (*$\k$*). (*$\c$*) ((*$\A\B$*) (*$\k$*)) ((*$\C\D$*) (*$\i$*))) (*$\dot$*) ((*$\lambda$*) (*$\k$*). (*$\c$*) (*$\B$*) ((*$\A\C\D$*) (*$\i$*) (*$\k$*))))
        (((*$\lambda$*) (*$\k$*). (*$\c$*) ((*$\A\C\D$*) (*$\i$*) (*$\k$*)) (*$\B$*)) (*$\dot$*) ((*$\lambda$*) (*$\k$*). (*$\c$*) ((*$\C\D$*) (*$\i$*)) ((*$\A\B$*) (*$\k$*)))))
    (((*$\lambda$*) (*$\j$*). ((*$\lambda$*) (*$\k$*). (*$\c$*) ((*$\A\B\C$*) (*$\j$*) (*$\k$*)) (*$\D$*)) (*$\dot$*) ((*$\lambda$*) (*$\k$*). (*$\c$*) ((*$\B\C$*) (*$\j$*)) ((*$\A\D$*) (*$\k$*))))
      (*$\dot$*) ((*$\lambda$*) (*$\j$*). ((*$\lambda$*) (*$\k$*). (*$\c$*) ((*$\A\C$*) (*$\k$*)) ((*$\B\D$*) (*$\j$*))) (*$\dot$*) ((*$\lambda$*) (*$\k$*). (*$\c$*) (*$\C$*) ((*$\A\B\D$*) (*$\j$*) (*$\k$*)))))
    (((*$\lambda$*) (*$\j$*). ((*$\lambda$*) (*$\k$*). (*$\c$*) ((*$\A\B\D$*) (*$\j$*) (*$\k$*)) (*$\C$*)) (*$\dot$*) ((*$\lambda$*) (*$\k$*). (*$\c$*) ((*$\B\D$*) (*$\j$*)) ((*$\A\C$*) (*$\k$*))))
      (*$\dot$*) ((*$\lambda$*) (*$\j$*). ((*$\lambda$*) (*$\k$*). (*$\c$*) ((*$\A\D$*) (*$\k$*)) ((*$\B\C$*) (*$\j$*))) (*$\dot$*) ((*$\lambda$*) (*$\k$*). (*$\c$*) (*$\D$*) ((*$\A\B\C$*) (*$\j$*) (*$\k$*)))))
\end{lstlisting}
Finally, \emph{just to contradict the epigraph}, we write the coherence for the $\cred{5}$-permutahedron:
\begin{lstlisting}
((*$\Perm$*) (*$\X$*) (*$\c$*))(*$_\cred{4}$*) (*$\A$*) (*$\B$*) (*$\A\B$*) (*$\C$*) (*$\A\C$*) (*$\B\C$*) (*$\A\B\C$*) (*$\D$*) (*$\A\D$*) (*$\B\D$*) (*$\A\B\D$*) (*$\C\D$*) (*$\A\C\D$*) (*$\B\C\D$*) (*$\A\B\C\D$*)
  (*$\E$*) (*$\A\E$*) (*$\B\E$*) (*$\A\B\E$*) (*$\C\E$*) (*$\A\C\E$*) (*$\B\C\E$*) (*$\A\B\C\E$*) (*$\D\E$*) (*$\A\D\E$*) (*$\B\D\E$*) (*$\A\B\D\E$*) (*$\C\D\E$*) (*$\A\C\D\E$*) (*$\B\C\D\E$*) :=
  (*$\Z\d\d\d\d$*) ((*$\S\d\d\d$*) ((*$\S\d\d$*) ((*$\S\d$*) ((*$\S$*) ((*$\Perm$*) (*$\X$*) (*$\c$*)) (*$\A$*)) (*$\B$*) (*$\A\B$*)) (*$\C$*) (*$\A\C$*) (*$\B\C$*) (*$\A\B\C$*))
    (*$\D$*) (*$\A\D$*) (*$\B\D$*) (*$\A\B\D$*) (*$\C\D$*) (*$\A\C\D$*) (*$\B\C\D$*) (*$\A\B\C\D$*)) (*$\E$*) (*$\A\E$*) (*$\B\E$*) (*$\A\B\E$*) (*$\C\E$*) (*$\A\C\E$*) (*$\B\C\E$*) (*$\A\B\C\E$*) (*$\D\E$*)
  (*$\A\D\E$*) (*$\B\D\E$*) (*$\A\B\D\E$*) (*$\C\D\E$*) (*$\A\C\D\E$*) (*$\B\C\D\E$*) =
  (*$\PathP$*)
    ((*$\lambda$*) (*$\i$*).
      (*$\PathP$*)
        ((*$\lambda$*) (*$\j$*).
          (*$\PathP$*)
            ((*$\lambda$*) (*$\k$*).
              (*$\Path$*) (*$\X$*) ((*$\c$*) (*$\A$*) ((*$\B\C\D\E$*) (*$\i$*) (*$\j$*) (*$\k$*))) ((*$\c$*) ((*$\B\C\D\E$*) (*$\i$*) (*$\j$*) (*$\k$*)) (*$\A$*)))
            (((*$\lambda$*) (*$\l$*). (*$\c$*) ((*$\A\B$*) (*$\l$*)) ((*$\C\D\E$*) (*$\i$*) (*$\j$*))) (*$\dot$*) ((*$\lambda$*) (*$\l$*). (*$\c$*) (*$\B$*) ((*$\A\C\D\E$*) (*$\i$*) (*$\j$*) (*$\l$*))))
            (((*$\lambda$*) (*$\l$*). (*$\c$*) ((*$\A\C\D\E$*) (*$\i$*) (*$\j$*) (*$\l$*)) (*$\B$*)) (*$\dot$*) ((*$\lambda$*) (*$\l$*). (*$\c$*) ((*$\C\D\E$*) (*$\i$*) (*$\j$*)) ((*$\A\B$*) (*$\l$*)))))
        (((*$\lambda$*) (*$\k$*). ((*$\lambda$*) (*$\l$*). (*$\c$*) ((*$\A\B\C$*) (*$\k$*) (*$\l$*)) ((*$\D\E$*) (*$\i$*))) (*$\dot$*) ((*$\lambda$*) (*$\l$*). (*$\c$*) ((*$\B\C$*) (*$\k$*)) ((*$\A\D\E$*) (*$\i$*) (*$\l$*))))
          (*$\dot$*) ((*$\lambda$*) (*$\k$*). ((*$\lambda$*) (*$\l$*). (*$\c$*) ((*$\A\C$*) (*$\l$*)) ((*$\B\D\E$*) (*$\i$*) (*$\k$*))) (*$\dot$*) ((*$\lambda$*) (*$\l$*). (*$\c$*) (*$\C$*) ((*$\A\B\D\E$*) (*$\i$*) (*$\k$*) (*$\l$*)))))
        (((*$\lambda$*) (*$\k$*). ((*$\lambda$*) (*$\l$*). (*$\c$*) ((*$\A\B\D\E$*) (*$\i$*) (*$\k$*) (*$\l$*)) (*$\C$*)) (*$\dot$*) ((*$\lambda$*) (*$\l$*). (*$\c$*) ((*$\B\D\E$*) (*$\i$*) (*$\k$*)) ((*$\A\C$*) (*$\l$*))))
          (*$\dot$*) ((*$\lambda$*) (*$\k$*). ((*$\lambda$*) (*$\l$*). (*$\c$*) ((*$\A\D\E$*) (*$\i$*) (*$\l$*)) ((*$\B\C$*) (*$\k$*))) (*$\dot$*) ((*$\lambda$*) (*$\l$*). (*$\c$*) ((*$\D\E$*) (*$\i$*)) ((*$\A\B\C$*) (*$\k$*) (*$\l$*))))))
    (((*$\lambda$*) (*$\j$*). ((*$\lambda$*) (*$\k$*). ((*$\lambda$*) (*$\l$*). (*$\c$*) ((*$\A\B\C\D$*) (*$\j$*) (*$\k$*) (*$\l$*)) (*$\E$*)) (*$\dot$*) ((*$\lambda$*) (*$\l$*). (*$\c$*) ((*$\B\C\D$*) (*$\j$*) (*$\k$*)) ((*$\A\E$*) (*$\l$*))))
        (*$\dot$*) ((*$\lambda$*) (*$\k$*). ((*$\lambda$*) (*$\l$*). (*$\c$*) ((*$\A\C\D$*) (*$\j$*) (*$\l$*)) ((*$\B\E$*) (*$\k$*))) (*$\dot$*) ((*$\lambda$*) (*$\l$*). (*$\c$*) ((*$\C\D$*) (*$\j$*)) ((*$\A\B\E$*) (*$\k$*) (*$\l$*)))))
      (*$\dot$*) ((*$\lambda$*) (*$\j$*). ((*$\lambda$*) (*$\k$*). ((*$\lambda$*) (*$\l$*). (*$\c$*) ((*$\A\B\D$*) (*$\k$*) (*$\l$*)) ((*$\C\E$*) (*$\j$*))) (*$\dot$*) ((*$\lambda$*) (*$\l$*). (*$\c$*) ((*$\B\D$*) (*$\k$*)) ((*$\A\C\E$*) (*$\j$*) (*$\l$*))))
        (*$\dot$*) ((*$\lambda$*) (*$\k$*). (((*$\lambda$*) (*$\l$*). (*$\c$*) ((*$\A\D$*) (*$\l$*)) ((*$\B\C\E$*) (*$\j$*) (*$\k$*))) (*$\dot$*) ((*$\lambda$*) (*$\l$*). (*$\c$*) (*$\D$*) ((*$\A\B\C\E$*) (*$\j$*) (*$\k$*) (*$\l$*)))))))
    (((*$\lambda$*) (*$\j$*). ((*$\lambda$*) (*$\k$*). ((*$\lambda$*) (*$\l$*). (*$\c$*) ((*$\A\B\C\E$*) (*$\j$*) (*$\k$*) (*$\l$*)) (*$\D$*)) (*$\dot$*) ((*$\lambda$*) (*$\l$*). (*$\c$*) ((*$\B\C\E$*) (*$\j$*) (*$\k$*)) ((*$\A\D$*) (*$\l$*))))
        (*$\dot$*) ((*$\lambda$*) (*$\k$*). ((*$\lambda$*) (*$\l$*). (*$\c$*) ((*$\A\C\E$*) (*$\j$*) (*$\l$*)) ((*$\B\D$*) (*$\k$*))) (*$\dot$*) ((*$\lambda$*) (*$\l$*). (*$\c$*) ((*$\C\E$*) (*$\j$*)) ((*$\A\B\D$*) (*$\k$*) (*$\l$*)))))
      (*$\dot$*) ((*$\lambda$*) (*$\j$*). ((*$\lambda$*) (*$\k$*). ((*$\lambda$*) (*$\l$*). (*$\c$*) ((*$\A\B\E$*) (*$\k$*) (*$\l$*)) ((*$\C\D$*) (*$\j$*))) (*$\dot$*) ((*$\lambda$*) (*$\l$*). (*$\c$*) ((*$\B\E$*) (*$\k$*)) ((*$\A\C\D$*) (*$\j$*) (*$\l$*))))
        (*$\dot$*) ((*$\lambda$*) (*$\k$*). (((*$\lambda$*) (*$\l$*). (*$\c$*) ((*$\A\E$*) (*$\l$*)) ((*$\B\C\D$*) (*$\j$*) (*$\k$*))) (*$\dot$*) ((*$\lambda$*) (*$\l$*). (*$\c$*) (*$\E$*) ((*$\A\B\C\D$*) (*$\j$*) (*$\k$*) (*$\l$*)))))))
\end{lstlisting}

\noindent
These $\n$-simplex types are obtained by directly applying the computation rules of Displayed Type Theory to the not-strictly-well-typed definition. This calculation has been implemented in a \texttt{Racket} program. The Appendix to this paper explains both the theory and practice of calculating $\n$-simplex types. There, we can find a link to a \texttt{Github} repository with the code, as well as the outputted coherences, up to the case of a $\cred{9}$-permutahedron.

Note that the code alone does not account for path lengths. Thus, to obtain an translation of the results in the syntax for cubical Moore paths introduced in the interlude, we simply need to add path length indices ranging from $\cred{1}$ to $\n$ from outermost to innermost $\PathP$/$\Path$.

Do these formulas become well-typed in the setting of cubical Moore paths?
Well, since $\bkt{\Perm\;\X\;\c}_\cred{2}$ being well-typed requires $\c$ to be definitionally associative. It stands to reason that $\bkt{\Perm\;\X\;\c}_\cred{3}$ being well-typed requires the dependent composition operation $\c\cblu{'}$ to be definitionally dependently associative. Looking at \ref{eq:dep-assoc}, we can write this out in the simpler case without explicit lower association operations:
\begin{align*}
\bkt{\c\cblu{'}\;\B\;\A\B\;\bkt{\c\;\C\;\D}\;\bkt{\c\cblu{'}\;\C\;\A\C\;\D\;\A\D}} \equiv \bkt{\c\cblu{'}\;\bkt{\c\cblu{'}\;\B\;\C}\;\bkt{\c\cblu{'}\;\B\;\A\B\;\C\;\A\C}\; \D\;\A\D}
\end{align*}
Which says that:
\begin{align*}
&\bkt{\lambda\;\i\text{.}~ \c\;\bkt{\A\B\;\i}\;\bkt{\c\;\C\;\D}} \dot \bkt{\lambda\;\i\text{.}~\c\;\B\;\bkt{\bkt{\bkt{\lambda\;\j\text{.}~ \c\;\bkt{\A\C\;\j}\;\D} \dot \bkt{\lambda\;\j\text{.}~\c\;\C\;\bkt{\A\D\;\j}}}\;\i}} \\
\equiv &\bkt{\lambda\;\i\text{.}~ \c\;\bkt{\bkt{\bkt{\lambda\;\j\text{.}~ \c\;\bkt{\A\B\;\j}\;\C} \dot \bkt{\lambda\;\j\text{.}~\c\;\B\;\bkt{\A\C\;\j}}}\;\i}\;\D} \dot \bkt{\lambda\;\i\text{.}~\c\;\bkt{\c\;\B\;\C}\;\bkt{\A\D\;\i}}
\end{align*}
In order for this to hold definitionally, in addition to the associativity of $\c$, we need path composition to be definitionally associative and for whiskering to distribute definitionally over path composition. More generally, $\bkt{\Perm\;\X\;\c}_\cred{4}$ being well-typed requires $\c\cblu{''}$ to be definitionally associative in the appropriate higher sense, and so on. All of these identities come down to path algebra that holds in the nice setting of cubical Moore paths.

We conclude by answering the problem posed at the end of the prelude: \emph{To define what it means for $\c$ to be commutative to homotopy in an infinitely coherent way.} Conceptually, this just means that we want a functorial assignment of an $\n$-simplex in $\Perm~\X~\c$ to each $\bkt{\n+\cred{1}}$ letters in $\X$. To state this formally, we may form a semi-simplicial type $\mathsf{\cred{coSk}}~\X : \SST$ from the type $\X : \Type$ by extending coskeletally. Then, a proof of the commutativity of $\c$ is a section of the projection $\cred{\mathsf{Hom}}\,\bkt{\Perm~\X~\c}\,\bkt{\mathsf{\cred{coSk}}~\X}$. \bbox

\begin{center}
\textbf{\cblu{\circled{\Large \cred{\stopbar}}}\ \ {\large Finale}}
\end{center}

\begin{center}
I'm not saying that \emph{You Can't Permutahedron}; I'm saying that \emph{You Wouldn't Permutahedron}.\\
And the only appropriate response to being told you that you wouldn't do something is:
\vspace{-0.6cm}

\[\textit{`Watch me!'} \tag*{\bbox}\] 
\end{center}
\vspace{-0.4cm}

\bibliographystyle{alpha}
\bibliography{main}

\clearpage

\begin{center}
\textbf{{\large Appendix: Calculating Simplex Types}}
\end{center}

In this appendix, we present a simplified version of \emph{Displayed Type Theory (dTT)}, where, in particular, we omit the modal aspects of the type theory and use Russell style universes. We are concerned not only with the theory of display and semi-simplicial types, but also with the practice of implementing a rather narrow \texttt{Racket} program for calculating the $\n$-simplex types of $\Perm~\X~\c$.
\vspace{0.2cm}

\noindent
\textbf{Introduction}

Consider the standard proof for the normalisation of Simply Typed Lambda Calculus in terms of computability predicates. For every type $\A$ of the theory, we define a predicate on the terms $\t : \A$, denoted $\dangle{\A}~\t$. At the base type, we have that $\dangle{\Base}~\t$ holds just in case $\t$ normalises. At arrow types, this predicate is defined recursively; we say that $\dangle{\A \to \B}~\f$ holds just in case $\f$ normalises, and if for every $\a : \A$ such that $\dangle{\A}~\a$ holds, we have that $\dangle{\B}~\bkt{\f~\a}$ holds. The proof then proceeds by structural induction on the syntax of $\t : \A$ to demonstrate that we always have that $\dangle{\A}~\t$ holds.

The phenomenon in this proof is known as \emph{unary parametricity}. In the above, we had that the predicates were metatheoretical as the statement \emph{`$\t$ normalises'} is not an internal construction to type theory. However, we may also consider theories that \emph{internalise} the parametricity operation. In the case of dTT, this internalisation takes the following form -- for every $\A : \Type$ and $\t : \A$, then we get $\A\d~\t : \Type$ and $\t\d : \A\d~\t$. That is, every type $\A : \Type$ defines a predicate $\A\d$ on its terms, and every term $\t : \A$ satisfies its corresponding predicate, such that we have a proof $\t\d : \A\d~\t$. The operation $\d$ is known as \emph{display}. In the case of arrow types, we have that $\bkt{\A\to\B}\d~\f \equiv \bkt{\a:\A} \to \A\d~\a \to \B\d~\bkt{\f~\a}$, which mirrors the recursive definition of the computability predicate above. \bbox
\vspace{0.2cm}

\noindent
\textbf{Meta-Abstraction}

Displayed Type Theory differed from all other previously existing internal parametricity theories in that its parametricity was only \emph{partially} internal. In particular, display is an example of a new notion called an \emph{indexed modality}. In general, a modality is a non-binding operation that changes the context. For example, in the case the formulation of dependent modal logic using locks, if we want to form the type $\sq~\A$ in context $\Gamma$, that is, if we want to make the judgement $\Gamma \vdash \sq~\A : \Type$, then we must have $\Gamma,~\lock_\sq \vdash \A : \Type$. Morally, the lock is a restriction on variable use rules in $\A$, and \emph{to first approximation} the lock modifies $\Gamma$ by getting rid of certain variables.

Display is a much more insidious kind of modality. If we have $\Gamma \vdash \A : \Type$, then $\Gamma\DD \vdash \A\d : \A \to \Type$. Here $\Gamma\DD$ is a context of double the length of $\Gamma$ in which every variable $\big(\a : \A\big)$ is replaced with $\bkt{\a : \A,~\a\cblu{'} : \A\d~\a}$, that is, the newly introduced variables are the parametricity witnesses of the existing variables. Thus, instead of \emph{deleting} variables from the context, as with the $\sq$ modality, display \emph{introduces} new variables without giving them a binding, which is very problematic. dTT deals with this issue by introducing modal guards, a theory which occupies many pages of the dTT paper.

From the perspective of a programmer, however, a much simpler version of the theory just states that display may only be applied to closed terms. Hence:
\begin{mathpar}
\infer{\vdash \A : \Type}{\Gamma \vdash \A\d : \A \to \Type}
\and
\infer{\vdash \t : \A}{\Gamma \vdash \t\d : \A\d~\t}
\end{mathpar}
In fact, we have $\Type\d \equiv \lambda~\A.~\A\to\Type$, so the first rule is a special case of the second.

However, from the perspective of the implementer of the programming language, we would like display to compute by structural recursion on the term, which would involve going under binders. Thus, although the user is not allowed to call display on open terms, the language implementation must.

In this appendix, we use Russel style universes and pun terms and types, which differs from usual theoretical practice. We will thus write the judgement $\Gamma \vdash \t \term$ to say that $\t$ is a term of some type (possibly of the universe). We then represent the notion of an open term with a \emph{meta-abstraction}:
\begin{mathpar}
\infer{\gamma : \Gamma \vdash \t \term}{\dsqbkt{\,\t\,}_{\,\gamma\,:\,\Gamma} \slfrac{\term}{_{\gamma\,:\,\Gamma}}}
\and
\infer{\Gamma \vdash \t : \A}{\dsqbkt{\,\t\,}_{\,\gamma\,:\,\Gamma} : \dsqbkt{\,\A\,}_{\,\gamma\,:\,\Gamma}}
\end{mathpar}
Meta-abstractions themselves are not first-class terms of the programming language, and must be fully applied to a list of terms to become one:
\begin{mathpar}
\infer{\cblu{\mathfrak{t}} \slfrac{\term}{_{\delta\,:\,\Delta}} \\ \sigma : \Gamma \Rightarrow \Delta}{\Gamma \vdash \cblu{\mathfrak{t}}~\sigma \term}
\and
\infer{\cblu{\mathcal{A}} \slfrac{\term}{_{\delta\,:\,\Delta}} \\ \cblu{\mathfrak{t}} : \cblu{\mathcal{A}} \\ \sigma : \Gamma \Rightarrow \Delta}{\Gamma \vdash \cblu{\mathfrak{t}}~\sigma : \cblu{\mathcal{A}}~\sigma}
\end{mathpar}
We may then identify closed terms $\t$ with their meta-abstractions $\dsqbkt{\,\t\,}_{\,\sqbkt{\,} : \bkt{}}$. \bbox
\vspace{0.2cm}

\noindent
\textbf{Display on Meta-Abstractions}

Two fundamental ingredients that we will need are \emph{d\'ecalage} and the \emph{evens projection}:
\begin{mathpar}
\infer{\Gamma \ctx}{\Gamma\DD \ctx}
\and
\infer{\sigma^\cblu{+} : \Gamma \Rightarrow \Delta\DD}{\sigma^\cblu{+}\ev : \Gamma \Rightarrow \Delta}
\end{mathpar}
Display on meta-abstractions then takes the following form:
\begin{mathpar}
\infer{\cblu{\mathfrak{t}} \slfrac{\term}{_{\gamma\,:\,\Gamma}}}{\cblu{\mathfrak{t}}\d \slfrac{\term}{_{\gamma^\cblu{+}\,:\,\Gamma\DD}}}
\end{mathpar}
The first rule that we will give for computing display is as follows:
\[\dsqbkt{\,\Type\,}_{\,\gamma\,:\,\Gamma}\d \equiv \dsqbkt{\,\lambda~\bkt{\A : \Type}\text{.}~\A \to \Type\,}_{\,\gamma^\cblu{+}\,:\,\Gamma\DD}\]
In terms of typing judgements, we then have:
\begin{mathpar}
\infer{\cblu{\mathfrak{t}} \slfrac{\term}{_{\gamma\,:\,\Gamma}} \\ \cblu{\mathcal{A}} \slfrac{\term}{_{\gamma\,:\,\Gamma}} \\ \cblu{\mathfrak{t}} : \cblu{\mathcal{A}}}{\cblu{\mathfrak{t}}\d : \dsqbkt{\,\cblu{\mathcal{A}}\d~\gamma^\cblu{+}~\bkt{\cblu{\mathfrak{t}}~\gamma^\cblu{+}\ev}\,}_{\,\gamma^\cblu{+}\,:\,\Gamma\DD}}
\end{mathpar}
Note that the application of $\cblu{\mathfrak{t}}~\gamma^\cblu{+}\ev$ makes sense, since $\cblu{\mathfrak{t}} : \cblu{\mathcal{A}}$ means $\cblu{\mathcal{A}} : \dsqbkt{\,\Type\,}_{\,\gamma\,:\,\Gamma}$, and thus that $\cblu{\mathcal{A}}\d : \dsqbkt{\,\cblu{\mathcal{A}}~\gamma^\cblu{+}\ev \to \Type\,}_{\,\gamma^\cblu{+}\,:\,\Gamma\DD}$, so after we instantiate $\cblu{\mathcal{A}}\d$ on $\gamma^\cblu{+}$, we may apply $\cblu{\mathfrak{t}}~\gamma^\cblu{+}\ev : \cblu{\mathcal{A}}~\gamma^\cblu{+}\ev$. From here, the definitions of d\'ecalage and the evens projection are as follows:
\begin{align*}
\big(\gamma:\Gamma,~\a : \A\big)\DD &\equiv \big(\gamma^\cblu{+} : \Gamma\DD, ~\a : \A~\sqbkt{\,\gamma^\cblu{+}\ev\,}, ~\a\cblu{'} : \dsqbkt{\,\A\,}_{\,\gamma\,:\,\Gamma}\d~\gamma^\cblu{+}~\a\big) \\
\sqbkt{\,\sigma^\cblu{+},~\t,~\t\cblu{'}\,}\ev &\equiv \sqbkt{\,\sigma^\cblu{+}\ev,~\t\,}
\end{align*}
For the rest of the term formers encountered in this paper, the computation rules for display are as follows:
\begin{align*}
\scaleto{\newllbracket}{1.2em}\,\bkt{\a : \A} \to \B\,\scaleto{\newrrbracket}{1.2em}_{\,\gamma\,:\,\Gamma}\d &\equiv \scaleto{\newllbracket}{1.2em}\, \lambda~\f\text{.}~\big(\a : \A~\sqbkt{\,\sigma^\cblu{+}\ev\,}\big)\,\big(\a\cblu{'} : \dsqbkt{\,\A\,}_{\,\gamma\,:\,\Gamma}\d~\gamma^\cblu{+}~\a\big) \to \\
&\quad\quad\quad\quad\dsqbkt{\,\B\,}_{\,\gamma\,:\,\Gamma,~\a\,:\,\A}\d~\gamma^\cblu{+}~\a~\a\cblu{'}~\bkt{\f~\a}\,\scaleto{\newrrbracket}{1.2em}_{\,\gamma^\cblu{+}\,:\,\Gamma\DD} \\
\scaleto{\newllbracket}{1.2em}\, \lambda~\a\text{.}~\t\,\scaleto{\newrrbracket}{1.2em}_{\,\gamma\,:\,\Gamma}\d &\equiv \scaleto{\newllbracket}{1.2em}\, \lambda~\a\text{.}~\lambda.~\a\cblu{'}\text{.}~\dsqbkt{\,\t\,}_{\,\gamma\,:\,\Gamma,~\a\,:\,\A}\d~\gamma^\cblu{+}~\a~\a\cblu{'} \,\scaleto{\newrrbracket}{1.2em}_{\,\gamma^\cblu{+}\,:\,\Gamma\DD} \\
\scaleto{\newllbracket}{1.2em}\, \f~\a \,\scaleto{\newrrbracket}{1.2em}_{\,\gamma\,:\,\Gamma}\d &\equiv \scaleto{\newllbracket}{1.2em}\, \big(\dsqbkt{\,\f\,}_{\,\gamma\,:\,\Gamma}\d~\gamma^\cblu{+}\big)~\big(\a~\sqbkt{\,\gamma^\cblu{+}\ev\,}\big)~\big(\dsqbkt{\,\a\,}_{\,\gamma\,:\,\Gamma}\d~\gamma^\cblu{+}\big) \,\scaleto{\newrrbracket}{1.2em}_{\,\gamma^\cblu{+}\,:\,\Gamma\DD} \\
\scaleto{\newllbracket}{1.2em}\, \Path~\X~\a~\b \,\scaleto{\newrrbracket}{1.2em}_{\,\gamma\,:\,\Gamma}\d &\equiv \scaleto{\newllbracket}{1.2em}\, \lambda~\alpha\text{.}~\PathP~\big(\lambda~\i\text{.}~\dsqbkt{\,\X\,}_{\,\gamma\,:\,\Gamma}\d~\gamma^\cblu{+} ~\bkt{\alpha~\i}\big) \\
&\quad\quad\quad\quad\big(\dsqbkt{\,\a\,}_{\,\gamma\,:\,\Gamma}\d~\gamma^\cblu{+}\big) ~\big(\dsqbkt{\,\b\,}_{\,\gamma\,:\,\Gamma}\d~\gamma^\cblu{+}\big)\,\scaleto{\newrrbracket}{1.2em}_{\,\gamma^\cblu{+}\,:\,\Gamma\DD} \\
\scaleto{\newllbracket}{1.2em}\, \PathP~\bkt{\lambda~\i\text{.}~\X}~\a~\b \,\scaleto{\newrrbracket}{1.2em}_{\,\gamma\,:\,\Gamma}\d &\equiv \scaleto{\newllbracket}{1.2em}\, \lambda~\alpha\text{.}~\PathP~\big(\lambda~\i\text{.}~\dsqbkt{\,\X\,}_{\,\gamma\,:\,\Gamma,~\i\,:\,\I}\d~\gamma^\cblu{+}~\i~\bkt{\alpha~\i}\big) \\
&\quad\quad\quad\quad\big(\dsqbkt{\,\a\,}_{\,\gamma\,:\,\Gamma}\d~\gamma^\cblu{+}\big) ~\big(\dsqbkt{\,\b\,}_{\,\gamma\,:\,\Gamma}\d~\gamma^\cblu{+}\big)\,\scaleto{\newrrbracket}{1.2em}_{\,\gamma^\cblu{+}\,:\,\Gamma\DD} \\
\scaleto{\newllbracket}{1.2em}\, \lambda~\i\text{.}~\t \,\scaleto{\newrrbracket}{1.2em}_{\,\gamma\,:\,\Gamma}\d &\equiv \scaleto{\newllbracket}{1.2em}\, \lambda~\i\text{.}~\dsqbkt{\,\t\,}_{\,\gamma\,:\,\Gamma,~\i\,:\,\I}\d~\gamma^\cblu{+}~\i \,\scaleto{\newrrbracket}{1.2em}_{\,\gamma^\cblu{+}\,:\,\Gamma\DD} \\
\scaleto{\newllbracket}{1.2em}\, \alpha~\i \,\scaleto{\newrrbracket}{1.2em}_{\,\gamma\,:\,\Gamma}\d &\equiv \scaleto{\newllbracket}{1.2em}\, \dsqbkt{\,\alpha\,}_{\,\gamma\,:\,\Gamma}\d~\gamma^\cblu{+}~\bkt{\i~\sqbkt{\,\gamma^\cblu{+}\ev\,}} \,\scaleto{\newrrbracket}{1.2em}_{\,\gamma^\cblu{+}\,:\,\Gamma\DD} \\
\scaleto{\newllbracket}{1.2em}\,\alpha~\dot~\beta \,\scaleto{\newrrbracket}{1.2em}_{\,\gamma\,:\,\Gamma}\d &\equiv \scaleto{\newllbracket}{1.2em}\, \big(\dsqbkt{\,\alpha\,}_{\,\gamma\,:\,\Gamma}\d~\gamma^\cblu{+}\big)  \dot \big(\dsqbkt{\,\beta\,}_{\,\gamma\,:\,\Gamma}\d~\gamma^\cblu{+}\big)\,\scaleto{\newrrbracket}{1.2em}_{\,\gamma^\cblu{+}\,:\,\Gamma\DD}
\end{align*}
These rules suggest that the cubical interval \emph{pretype} $\I$ behaves differently from a regular type in terms of display. The general rule is that we do not allow forming the displayed type $\I\d$, and, in the course of d\'ecalage, interval variables do not get doubled, with the definition of the evens projection being adjusted accordingly. We see that we must also distinguish abstractions and applications of interval variables from regular abstractions and applications. Note that the rule for display on $\PathP$ types assumes that the line of types is given by an interval abstraction, which may generally be assumed by way of $\eta$-expansion.

The final rule to give concerns display of a variable:
\[\scaleto{\newllbracket}{1.2em}\, \x \,\scaleto{\newrrbracket}{1.2em}_{\,\gamma\,:\,\Gamma}\d \equiv \scaleto{\newllbracket}{1.2em}\, \x\cblu{'}\,\scaleto{\newrrbracket}{1.2em}_{\,\gamma^\cblu{+}\,:\,\Gamma\DD},\]
where $\x\cblu{'}$ is the twinned variable arising from d\'ecalage. In the case of interval variables, we \emph{could} say there is no twinned variable and display is inert:
\[\scaleto{\newllbracket}{1.2em}\, \i \,\scaleto{\newrrbracket}{1.2em}_{\,\gamma\,:\,\Gamma}\d \equiv \scaleto{\newllbracket}{1.2em}\, \i~\sqbkt{\,\gamma^\cblu{+}\ev\,} \,\scaleto{\newrrbracket}{1.2em}_{\,\gamma^\cblu{+}\,:\,\Gamma\DD},\]
however it would be \emph{better} to say that forming display of an interval variable is illegal. \bbox
\vspace{0.1cm}

\noindent
\textbf{Telescopes and Partial Substitutions}

As a minor extension of the above theory, we may consider the notions of \emph{telescopes} and \emph{partial substitutions}:
\begin{mathpar}
\infer{\bkt{\Gamma,~\Delta} \ctx}{\Gamma \vdash \Delta \tel}
\and
\infer{\sqbkt{\,\cred{1}_\Gamma,~\sigma\,} : \Gamma \Rightarrow \bkt{\Gamma,~\Delta}}{\Gamma \vdash \sigma : \Delta}
\end{mathpar}
These may then be meta-abstracted to form:
\begin{mathpar}
\infer{\Gamma \vdash \Delta \tel}{\dsqbkt{\,\Delta\,}_{\,\gamma\,:\,\Gamma} \slfrac{\tel}{_{\gamma\,:\,\Gamma}}}
\and
\infer{\Gamma \vdash \sigma : \Delta}{\dsqbkt{\,\sigma\,}_{\,\gamma\,:\,\Gamma} : \dsqbkt{\,\Delta\,}_{\,\gamma\,:\,\Gamma}}
\end{mathpar}
And the operation with which we are concerned is d\'ecalage:
\begin{mathpar}
\infer{\mathcal{\cblu{D}} \slfrac{\tel}{_{\gamma\,:\,\Gamma}}}{\mathcal{\cblu{D}}\DD \slfrac{\tel}{_{\gamma^\cblu{+}\,:\,\Gamma\DD}}}
\and
\infer{\mathcal{\cblu{D}} \slfrac{\tel}{_{\gamma\,:\,\Gamma}} \\ \cblu{\mathfrak{s}} : \mathcal{\cblu{D}}}{\cblu{\mathfrak{s}}\DD : \mathcal{\cblu{D}}\DD}
\end{mathpar}
For which the computation rules are:
\begin{align*}
\scaleto{\newllbracket}{1.2em} \bkt{} \scaleto{\newrrbracket}{1.2em}_{\,\gamma\,:\,\Gamma}\DD &\equiv \scaleto{\newllbracket}{1.2em} \bkt{} \scaleto{\newrrbracket}{1.2em}_{\,\gamma^\cblu{+}\,:\,\Gamma\DD} \\
\scaleto{\newllbracket}{1.2em} \big(\a : \A,~\delta : \Delta\big) \scaleto{\newrrbracket}{1.2em}_{\,\gamma\,:\,\Gamma}\DD &\equiv \scaleto{\newllbracket}{1.2em} \big(\a : \A~\sqbkt{\,\gamma^\cblu{+}\ev\,},~\a\cblu{'} : \dsqbkt{\,\A\,}_{\,\gamma\,:\,\Gamma}\d ~\gamma^\cblu{+}~\a, \\
&\quad\quad\quad\quad \delta^\cblu{+} : \dsqbkt{\,\Delta\,}_{\,\gamma\,:\,\Gamma,~\a\,:\,\A}\DD~\gamma^\cblu{+}~\a~\a\cblu{'}\big) \scaleto{\newrrbracket}{1.2em}_{\,\gamma^\cblu{+}\,:\,\Gamma\DD} \\
\scaleto{\newllbracket}{1.2em}\, \sqbkt{\,} \,\scaleto{\newrrbracket}{1.2em}_{\,\gamma\,:\,\Gamma}\DD &\equiv \scaleto{\newllbracket}{1.2em}\, \sqbkt{\,} \,\scaleto{\newrrbracket}{1.2em}_{\,\gamma^\cblu{+}\,:\,\Gamma\DD} \\
\scaleto{\newllbracket}{1.2em}\, \big[\,\t,~\sigma\,\big] \,\scaleto{\newrrbracket}{1.2em}_{\,\gamma\,:\,\Gamma}\DD &\equiv \scaleto{\newllbracket}{1.2em}\, \big[\,\t~\sqbkt{\,\gamma^\cblu{+}\ev\,},~\dsqbkt{\,\t\,}_{\,\gamma\,:\,\Gamma}\d~\gamma^\cblu{+},~\dsqbkt{\,\sigma\,}_{\,\gamma\,:\,\Gamma}\DD~\gamma^\cblu{+}\,\big] \,\scaleto{\newrrbracket}{1.2em}_{\,\gamma^\cblu{+}\,:\,\Gamma\DD}
\end{align*}
Note that telescopes are more naturally thought of as being constructed by iteratively bringing types over from $\Gamma$ to $\Delta$ in $\bkt{\Gamma,~\Delta}$ and thus as being constructed from the left, as opposed to from the right (as contexts are). \bbox
\vspace{0.2cm}

\noindent
\textbf{Semi-Simplicial Types}

The definition of semi-simplicial types in dTT is as follows:
\begin{lstlisting}
(*\textbf{codata}*) (*$\SST$*) : (*$\Type$*) (*\textbf{where}*)
  (*$\Z$*) : (*$\SST$*) (*$\to$*) (*$\Type$*)
  (*$\S$*) : ((*$\A$*) : (*$\SST$*)) (*$\to$*) (*$\Z$*) (*$\A$*) (*$\to$*) (*$\SST\d$*) (*$\A$*)
\end{lstlisting}
In the general theory, this is an instance of a larger class of definitions known as \emph{displayed coinductive types}, but we will focus only on the coinduction principle of $\SST$. The basic rules regarding the constants $\SST$, $\Z$, and $\S$ are:
\[\infer{\hspace{0.1cm}}{\Shortstack{{$\Gamma \vdash \SST : \Type$} \hspace{0.1cm} {$\Z : \dsqbkt{\,\Type\,}_{\,\A\;:\;\SST}$} \hspace{0.1cm} {$\S: \dsqbkt{\,\SST^{\mathsf{\cred{d}}}~\A\,}_{\,\A\;:\;\SST,\;\a\;:\Z\;\A}$}}}\]

If we consider reasonable looking code for a definition of an $\SST$ from a context of two parameters, it would look as follows:
\begin{lstlisting}
(*$\f$*) : ((*$\t$*) : (*$\X$*)) ((*$\s$*) : (*$\Y$*) (*$\t$*)) (*$\to$*) (*$\SST$*)
(*$\Z$*) ((*$\f$*) (*$\t$*) (*$\s$*)) = ((*\cprime{?$_{\mathsf{Z}_1}$}*) : (*$\Type$*))
(*$\S$*) ((*$\f$*) (*$\t$*) (*$\s$*)) (*$\a$*) = (*$\f\d$*) (*$\t$*) ((*\cprime{?$_{\mathsf{S}_1}$}*) : (*$\X\d$*) (*$\t$*)) (*$\s$*) ((*\cprime{?$_{\mathsf{S}_2}$}*) : (*$\Y\d$*) (*$\t$*) (*\cprime{?$_{\mathsf{S}_1}$}*) (*$\s$*))
\end{lstlisting}
Here, suppose that $\Gamma \ctx$. If we think of $\Gamma$ as a \emph{state space} and $\sqbkt{\,\sigma : \Gamma\,}$ as a \emph{state}. Then the above definition suggests that we are able to define $\f : \bkt{\gamma : \Gamma} \to \SST$ provided that we are able to provide two ingredients. First, we need a way of extracting $\sqbkt{\,\cblu{\widetilde{\mathsf{Z}}}~\sigma : \Type\,}$, a type of $\cred{0}$-simplices, from a state $\sigma$. Second, we need a way of extracting $\sqbkt{\,\cblu{\widetilde{\mathsf{S}}}~\sigma~\a : \Gamma\d~\sigma\,}$, a dependent section of $\Gamma$ over $\sigma$, from a state $\sigma$ and a $\cred{0}$-simplex $\sqbkt{\,\a : \cblu{\widetilde{\mathsf{Z}}}~\sigma\,}$. This suggests that a reasonable coinduction principle for $\SST$ is the following:
\begin{mathpar}
\infer{\Gamma \ctx \\ \cblu{\widetilde{\mathsf{Z}}} : \dsqbkt{\,\Type\,}_{\,\gamma\;:\;\Gamma} \\ \cblu{\widetilde{\mathsf{S}}} :  \dsqbkt{\,\Gamma\d~\gamma\,}_{\;\gamma\;:\;\Gamma,\;\a\;:\;\cblu{\widetilde{\mathsf{Z}}}\;\gamma}}{\mathsf{\R}_\Gamma~\cblu{\widetilde{\mathsf{Z}}} ~\cblu{\widetilde{\mathsf{S}}} : \dsqbkt{\,\SST\,}_{\,\gamma\;:\;\Gamma}}
\end{mathpar}
and that its computation rules should be:
\begin{align*}
\Z~\big((\R_\Gamma~\cblu{\widetilde{\mathsf{Z}}}~\cblu{\widetilde{\mathsf{S}}})~\sigma\big) &\equiv \cblu{\widetilde{\mathsf{Z}}}~\sigma \\
\S~\big((\R_\Gamma~\cblu{\widetilde{\mathsf{Z}}} ~\cblu{\widetilde{\mathsf{S}}})~\sigma\big)~\a &\equiv \big(\R_\Gamma~\cblu{\widetilde{\mathsf{Z}}} ~\cblu{\widetilde{\mathsf{S}}}\big)\d~\pair{\sigma}{\cblu{\widetilde{\mathsf{S}}}~\sigma~\a}.
\end{align*}
Where $\pair{\blank}{\blank}$ is an operation that intertwines terms in the section $\sqbkt{\,\sigma : \Gamma\,}$ with their parametricity witnesses $\sqbkt{\,\cblu{\widetilde{\mathsf{S}}}~\sigma~\a : \Gamma\d~\sigma\,}$ such as to produce $\sqbkt{\,\pair{\sigma}{\cblu{\widetilde{\mathsf{S}}}~\sigma~\a} : \Gamma\DD\,}$.

The reasonable question to ask is whether or not the expression $(\R_\Gamma~\cblu{\widetilde{\mathsf{Z}}} ~\cblu{\widetilde{\mathsf{S}}})\d$ appearing above can be computed. The answer is \emph{yes} in an extension of dTT with \emph{symmetries}, and \emph{no} in base dTT. However, in base dTT, we are able to say how $\Z\d$ and $\S\d$ (and higher analogues thereof) compute on this expression, which is sufficient to define all the computational behaviour necessary to extract $\n$-simplex types. In particular, if we consider the program:
\begin{lstlisting}
(*$\f$*) : ((*$\t$*) : (*$\X$*)) ((*$\s$*) : (*$\Y$*) (*$\t$*)) (*$\to$*) (*$\SST$*)
(*$\Z$*) ((*$\f$*) (*$\t$*) (*$\s$*)) = (*$\cblu{\mathfrak{z}}$*) (*$\t$*) (*$\s$*)
(*$\S$*) ((*$\f$*) (*$\t$*) (*$\s$*)) (*$\a$*) = (*$\f\d$*) (*$\t$*) ((*$\cblu{\mathfrak{s}_0}$*) (*$\t$*) (*$\s$*) (*$\a$*)) (*$\s$*) ((*$\cblu{\mathfrak{s}_1}$*) (*$\t$*) (*$\s$*) (*$\a$*))
\end{lstlisting}
Then we can compute display on each line of the definition to obtain:
\begin{lstlisting}
(*$\f\d$*) : ((*$\t$*) : (*$\X$*)) ((*$\t\cblu{'}$*) : (*$\X\d$*) (*$\t$*)) ((*$\s$*) : (*$\Y$*) (*$\t$*)) ((*$\s\cblu{'}$*) : (*$\Y\d$*) (*$\t$*) (*$\t\cblu{'}$*) (*$\s$*)) (*$\to$*) (*$\SST\d$*) ((*$\f$*) (*$\t$*) (*$\s$*))
(*$\Z\d$*) ((*$\f$*) (*$\t$*) (*$\t\cblu{'}$*) (*$\s$*) (*$\s\cblu{'}$*)) = (*$\lambda$*) (*$\a$*). (*$\cblu{\mathfrak{z}}\d$*) (*$\t$*) (*$\t\cblu{'}$*) (*$\s$*) (*$\s\cblu{'}$*) (*$\a$*)
(*$\S\d$*) ((*$\f$*) (*$\t$*) (*$\t\cblu{'}$*) (*$\s$*) (*$\s\cblu{'}$*)) (*$\a$*) (*$\a\cblu{'}$*) =
  (*$\f\d\d$*) (*$\t$*) (*$\t\cblu{'}$*) ((*$\cblu{\mathfrak{s}_0}$*) (*$\t$*) (*$\s$*) (*$\a$*)) ((*$\cblu{\mathfrak{s}_0}\d$*) (*$\t$*) (*$\t\cblu{'}$*) (*$\s$*) (*$\s\cblu{'}$*) (*$\a$*) (*$\a\cblu{'}$*)) (*$\s$*) (*$\s\cblu{'}$*) ((*$\cblu{\mathfrak{s}_1}$*) (*$\t$*) (*$\s$*) (*$\a$*)) ((*$\cblu{\mathfrak{s}_1}\d$*) (*$\t$*) (*$\t\cblu{'}$*) (*$\s$*) (*$\s\cblu{'}$*) (*$\a$*) (*$\a\cblu{'}$*))
\end{lstlisting}
Where we have chosen to suggestively $\eta$-expand the RHS of the second line. What this corresponds to is the following general rules:
\begin{align*}
\Z\d~\big((\R_\Gamma~\cblu{\widetilde{\mathsf{Z}}}~\cblu{\widetilde{\mathsf{S}}})\d~\sigma^\cblu{+}\big)~\a &\equiv \cblu{\widetilde{\mathsf{Z}}}\d~\sigma^\cblu{+}~\a \\
\S\d~\big((\R_\Gamma~\cblu{\widetilde{\mathsf{Z}}} ~\cblu{\widetilde{\mathsf{S}}})\d~\sigma^\cblu{+}\big)~\a~\a\cblu{'} &\equiv (\R_\Gamma~\cblu{\widetilde{\mathsf{Z}}} ~\cblu{\widetilde{\mathsf{S}}})\d\d~\big(\pair{\blank}{\blank}\DD~\sigma^\cblu{+}~(\cblu{\widetilde{\mathsf{S}}}\DD~\sigma^\cblu{+}~\a~\a\cblu{'})\big)
\end{align*}
Where the $\pair{\blank}{\blank}\DD$ operation concatenates groups of variables for each original variable in $\Gamma$, as above. Iterating this in general produces:
\begin{align*}
\Z^{\cred{\mathsf{d}}^\n}~\big((\R_\Gamma~\cblu{\widetilde{\mathsf{Z}}}~\cblu{\widetilde{\mathsf{S}}})^{\cred{\mathsf{d}}^\n}~\sigma^\n\big)~\cblu{\partial}\a &\equiv \cblu{\widetilde{\mathsf{Z}}}^{\cred{\mathsf{d}}^\n}~\sigma^\n~\cblu{\partial}\a \\
\S^{\cred{\mathsf{d}}^\n}~\big((\R_\Gamma~\cblu{\widetilde{\mathsf{Z}}} ~\cblu{\widetilde{\mathsf{S}}})^{\cred{\mathsf{d}}^\n}~\sigma^\n\big)~\cblu{\partial}\a~\a &\equiv (\R_\Gamma~\cblu{\widetilde{\mathsf{Z}}} ~\cblu{\widetilde{\mathsf{S}}})^{\cred{\mathsf{d}}^{\n+\cred{`}}}~\big(\pair{\blank}{\blank}^{\cred{\mathsf{D}}^\n}~\sigma^\n~(\cblu{\widetilde{\mathsf{S}}}^{\cred{\mathsf{D}}^\n}~\sigma^\n~\cblu{\partial}\a~\a)\big)
\end{align*}
This gives all of the rules necessary to compute the $\n$-simplex types of $\Perm~\X~\c$ and we will now make the theory much more concretely by way of implementation. \bbox
\vspace{0.2cm}

\noindent
\textbf{Implementation}

In the case of $\Perm$, we have:
\begin{align*}
\cred{\oldGamma} &\equiv \bkt{\X : \Type,~\c : \X \to \X \to \X} \\
\cred{\widetilde{\mathsf{Z}}} &\equiv \dsqbkt{\,\X\,}_{\,\X\,:\,\Type,~\c\,:\,\X\,\to\,\X\,\to\,\X} \\
\cred{\widetilde{\mathsf{S}}} &\equiv \scaleto{\newllbracket}{1.2em}\,\big[\,\lambda~\B\text{.}~\Path~\X~\bkt{\c~\A~\B}~\bkt{\c~\B~\A}, \\
&\quad\quad\ \ \lambda~\B~\A\B~\C~\A\C\text{.} \bkt{\lambda\;\i\text{.}~\c\;\bkt{\A\B\;\i}\;\C} \dot \bkt{\lambda\;\i\text{.}~\c\;\B\;\bkt{\A\C\;\i}}\,\big]\,\scaleto{\newrrbracket}{1.2em}_{\,\X\,:\,\Type,~\c\,:\,\X\,\to\,\X\,\to\,\X,~\A\,:\,\X}
\end{align*}
We will modify this slightly by splitting $\cred{\oldGamma} \equiv \bkt{\cred{\oldGamma_{X}}, ~\cred{\oldGamma_{\mu}}}$ and $\cred{\widetilde{\mathsf{S}}} \equiv \sqbkt{\,\cred{\widetilde{\S}_X},~\cred{\widetilde{\S}_\mu}\,}$, such that:
\begin{align*}
\cred{\oldGamma_{X}} &\equiv \dsqbkt{\,\X : \Type\,}_{\,\sqbkt{\,}\,:\,\bkt{}} \\
\cred{\oldGamma_{\mu}} &\equiv \dsqbkt{\,\c : \X \to \X \to \X\,}_{\,\X\,:\,\Type} \\
\cred{\widetilde{\S}_X} &\equiv \scaleto{\newllbracket}{1.2em}\,\big[\,\lambda~\B\text{.}~\Path~\X~\bkt{\c~\A~\B}~\bkt{\c~\B~\A}\,\big]\,\scaleto{\newrrbracket}{1.2em}_{\,\X\,:\,\Type,~\c\,:\,\X\,\to\,\X\,\to\,\X,~\A\,:\,\X} \\
\cred{\widetilde{\S}_\mu} &\equiv \scaleto{\newllbracket}{1.2em}\,\big[\,\lambda~\B~\A\B~\C~\A\C\text{.} \bkt{\lambda\;\i\text{.}~\c\;\bkt{\A\B\;\i}\;\C} \dot \bkt{\lambda\;\i\text{.} ~\c\;\B\;\bkt{\A\C\;\i}}\,\big]\,\scaleto{\newrrbracket}{1.2em}_{\,\X\,:\,\Type,~\c\,:\,\X\,\to\,\X\,\to\,\X,~\A\,:\,\X}
\end{align*}
Recall that we are trying to compute something like:
\[\Z\d\d\d~\big(\S\d\d~\big(\S\d~\big(\S~\bkt{\Perm~\X~\c}~\A\big)~\B~\A\B\big)~\C~\A\C~\B\C~\A\B\C\big)~\D~\A\D~\B\D~\A\B\D~\C\D~\A\C\D~\B\C\D\]
We thus start with the data of $\cred{\oldsigma_{X}^{0}} \equiv \sqbkt{\,\X\,} : \cred{\oldGamma_X}$ and $\cred{\oldsigma_{\mu}^{0}} \equiv \sqbkt{\,\c\,} : \cred{\oldGamma_\mu}\,\cred{\oldsigma_{X}^{0}}$, as well as several batches of letters $\sqbkt{\,\cred{\partial a}_\k,~\cred{a}_\k\,}$ where $\sqbkt{\,\cred{\partial a_0},~\cred{a_0}\,} \equiv \sqbkt{\,\A\,}$, $\sqbkt{\,\cred{\partial a_1},~\cred{a_1}\,} \equiv \sqbkt{\,\B,~\A\B\,}$, and $\sqbkt{\,\cred{\partial a_2},~\cred{a_2}\,} \equiv \sqbkt{\,\C,~\A\C,~\B\C,~\A\B\C\,}$, etc. In the case of the very last batch, we are only given $\sqbkt{\,\cred{\partial a}_\n\,}$, such that the last top dimensional group of letters is missing.

The gist of the computation is as follows: To extract the $\n$-simplex type $\bkt{\Perm~\X~\c}_\n$, we start with $\cred{\oldsigma_{X}^{0}}$ and $\cred{\oldsigma_{\mu}^{0}}$, and aim to iteratively compute $\cred{\oldsigma_{\cred{X}}^\n} : \cred{\oldGamma_X}^{\cred{\mathsf{D}}^\n}$ and $\cred{\oldsigma_{\cred{\mu}}^\n} : \cred{\oldGamma_\mu}^{\cred{\mathsf{D}}^\n}~\cred{\oldsigma_{X}^\n}$. We compute $\cred{\oldsigma}_{\cred{X}}^{\k+\cred{1}}$ and $\cred{\oldsigma}_{\cred{\mu}}^{\k+\cred{1}}$ as follows:
\begin{align*}
\cred{\oldsigma}_{\cred{X}}^{\k+\cred{1}} &\equiv \sqbkt{\,\cred{\oldsigma}_{\cred{X}}^\k,~\cred{\widetilde{\mathsf{S}}}_{\cred{X}}^{\cred{\mathsf{D}}^\k}~\cred{\oldsigma}_{\cred{X}}^\k~\cred{\oldsigma}_{\cred{\mu}}^\k~\cred{\partial a}_\k~\cred{a}_\k\,} \\
\cred{\oldsigma}_{\cred{\mu}}^{\k+\cred{1}} &\equiv \sqbkt{\,\cred{\oldsigma}_{\cred{\mu}}^\k,~\cred{\widetilde{\mathsf{S}}}_{\cred{\mu}}^{\cred{\mathsf{D}}^\k}~\cred{\oldsigma}_{\cred{X}}^\k~\cred{\oldsigma}_{\cred{\mu}}^\k~\cred{\partial a}_\k~\cred{a}_\k\,}.
\end{align*}
This computation corresponds to the iterative slicing. Now, since $\cred{\widetilde{\mathsf{Z}}}$ is just a meta-abstracted variable, we obtain $\Z^{\cred{\mathsf{d}}^\n}~\big((\Perm~\X~\c)^{\cred{\mathsf{d}}^\n}~\cred{\oldsigma}_{\cred{X}}^\n~\cred{\oldsigma}_{\cred{\mu}}^\n\big)$ simply by extracting the last term of $\cred{\oldsigma}_{\cred{X}}^\n$. Finally, we iteratively apply the variables in $\cred{\partial a}_\n$ to this term and normalise to obtain a formula for the desired $\n$-simplex type.

In terms of a \texttt{Racket} implementation, we will have five distinct kinds of variables:
\begin{lstlisting}
((*\cred{\textbf{struct}}*) (*\cred{XVar}*) ((*\cblu{level}*)))
((*\cred{\textbf{struct}}*) (*\cred{$\mu$Var}*) ((*\cblu{level}*)))
((*\cred{\textbf{struct}}*) (*\cred{$\oldalpha$Var}*) ((*\cblu{level}*)))
((*\cred{\textbf{struct}}*) (*\cred{$\oldGamma$Var}*) ((*\cblu{index}*)))
((*\cred{\textbf{struct}}*) (*\cred{iVar}*) ((*\cblu{index}*)))
\end{lstlisting}
The first three variable kinds represent the $\X$, $\c$, and $\A$ variable groups in an iterated d\'ecalage of the context $\bkt{\X : \Type,~\c : \X \to \X \to \X,~\A : \X}$, which is what will appear as a meta-binder in $\cred{\widetilde{\mathsf{S}}}_{\cred{X}}^{\cred{\mathsf{D}}^\k}$ and $\cred{\widetilde{\mathsf{S}}}_{\cred{\mu}}^{\cred{\mathsf{D}}^\k}$. For these we will use de Bruijn \emph{levels}. The last two kinds of variables represent the regular and interval variables acquired by going under concrete binders in a term.  For these we will use de Bruijn \emph{indices} and keep the regular and interval contexts separate.

We then have structures for all of the relevant term formers:
\begin{lstlisting}
((*\cred{\textbf{struct}}*) (*\cred{Lam}*) ((*\cblu{body}*)))
((*\cred{\textbf{struct}}*) (*\cred{Line}*) ((*\cblu{body}*)))
((*\cred{\textbf{struct}}*) (*\cred{App}*) ((*\cblu{fun}*) (*\cblu{arg}*)))
((*\cred{\textbf{struct}}*) (*\cred{iApp}*) ((*\cblu{fun}*) (*\cblu{arg}*)))
((*\cred{\textbf{struct}}*) (*\cred{Path}*) ((*\cblu{space}*) (*\cblu{lhs}*) (*\cblu{rhs}*)))
((*\cred{\textbf{struct}}*) (*\cred{PathP}*) ((*\cblu{line}*) (*\cblu{lhs}*) (*\cblu{rhs}*)))
((*\cred{\textbf{struct}}*) (*\cred{Comp}*) ((*\cblu{path1}*) (*\cblu{path2}*)))
\end{lstlisting}
Note that we do not have a structure for the type former of $\Pi$-types or the universe, since these type formers do not appear in our formulas.

Our first goal is to define renaming, substitution, and normalisation, generally following the approach in~\cite{LACI}. All of these will be 
roughly defined by structural recursion (although normalisation must make a generative call), and we will omit the uninteresting recursive calls. The latter two of these functions will only be called on terms, as opposed to meta-abstractions, and thus do not need to handle the case of the first three variables.
\clearpage

First, we cover the case of \emph{renaming}:
\begin{lstlisting}
((*\cred{\textbf{define}}*) ((*\cred{weaken-ren}*) (*\cblu{$\rho$}*))
  ((*\cred{$\lambda$}*) ((*\cblu{$n$}*)) ((*\cred{if}*) ((*\cred{zero?}*) (*\cblu{$n$}*)) (*\cred{$0$}*) ((*\cred{add1}*) ((*\cblu{$\rho$}*) ((*\cred{sub1}*) (*\cred{$n$}*)))))))

((*\cred{\textbf{define}}*) ((*\cred{rename}*) (*\cblu{expr}*) (*\cblu{$\rho$C}*) (*\cblu{$\rho$I}*))
  ((*\cred{\textbf{match}}*) (*\cblu{expr}*)
    [((*\cred{?}*) (*\cred{symbol?}*)) (*\cblu{expr}*)]
    [((*\cred{XVar}*) (*$\n$*)) ((*\cred{XVar}*) (*$\n$*))]
    [((*\cred{$\mu$Var}*) (*$\n$*)) ((*\cred{$\mu$Var}*) (*$\n$*))]
    [((*\cred{$\oldalpha$Var}*) (*$\n$*)) ((*\cred{$\oldalpha$Var}*) (*$\n$*))]
    [((*\cred{$\oldGamma$Var}*) (*$\n$*)) ((*\cred{$\oldGamma$Var}*) ((*\cblu{$\rho$C}*) (*$\n$*)))]
    [((*\cred{iVar}*) (*$\n$*)) ((*\cred{iVar}*) ((*\cblu{$\rho$I}*) (*$\n$*)))]
    [((*\cred{Lam}*) (*$\t$*)) ((*\cred{Lam}*) ((*\cred{rename}*) (*$\t$*) ((*\cred{weaken-ren}*) (*\cblu{$\rho$C}*)) (*\cblu{$\rho$I}*)))]
    [((*\cred{Line}*) (*$\t$*)) ((*\cred{Line}*) ((*\cred{rename}*) (*$\t$*) (*\cblu{$\rho$C}*) ((*\cred{weaken-ren}*) (*\cblu{$\rho$I}*))))]
    (*\cred{$\cdots$}*)))
\end{lstlisting}
The input to \cred{rename} is two functions $\cblu{\rho C}~\cblu{\rho I} : \N \to \N$, called \emph{renamings}, which dictate how the indices in the \cred{$\oldGamma$Var} and \cred{iVar} structure are to be transformed. Since we are using de Bruijn indices, we have to weaken when going under a binder. For example, if a renaming started off as $\bkt{\cred{2},\,\cred{5},\,\cred{7},\,\cdots}$, then weakening acts by shifting the sequence by appending a $\cred{0}$ and incrementing all subsequent entries, such as to obtain $\bkt{\cred{0},\,\cred{3},\,\cred{6},\,\cred{8},\,\cdots}$.

Second, we cover the case of \emph{substitution}:
\begin{lstlisting}
((*\cred{\textbf{define}}*) ((*\cred{weaken-sub}*) (*\cblu{$\rho$}*))
  ((*\cred{$\lambda$}*) ((*\cblu{$n$}*)) ((*\cred{if}*) ((*\cred{zero?}*) (*\cblu{$n$}*)) ((*\cred{$\oldGamma$Var}*) (*\cred{$0$}*)) ((*\cred{rename}*) ((*\cblu{$\rho$}*) ((*\cred{sub1}*) (*\cblu{$n$}*))) (*\cred{add1}*) ((*\cred{$\lambda$}*) ((*$\x$*)) (*$\x$*))))))

((*\cred{\textbf{define}}*) ((*\cred{weaken-sub-i}*) (*\cblu{$\rho$}*))
  ((*\cred{$\lambda$}*) ((*\cblu{$n$}*)) ((*\cred{rename}*) ((*\cblu{$\rho$}*) (*\cblu{$n$}*)) ((*\cred{$\lambda$}*) ((*$\x$*)) (*$\x$*)) (*\cred{add1}*))))

((*\cred{\textbf{define}}*) ((*\cred{subst}*) (*\cblu{expr}*) (*\cblu{$\rho$}*))
  ((*\cred{\textbf{match}}*) (*\cblu{expr}*)
    [((*\cred{?}*) (*\cred{symbol?}*)) (*\cblu{expr}*)]
    [((*\cred{$\oldGamma$Var}*) (*\cblu{$n$}*)) ((*\cblu{$\rho$}*) (*\cblu{$n$}*))]
    [((*\cred{iVar}*) (*\cblu{$n$}*)) ((*\cred{iVar}*) (*\cblu{$n$}*))]
    [((*\cred{Lam}*) (*$\t$*)) ((*\cred{Lam}*) ((*\cred{subst}*) (*$\t$*) ((*\cred{weaken-sub}*) (*\cblu{$\rho$}*))))]
    [((*\cred{Line}*) (*$\t$*)) ((*\cred{Line}*) ((*\cred{subst}*) (*$\t$*) ((*\cred{weaken-sub-i}*) (*\cblu{$\rho$}*))))]
    (*\cred{$\cdots$}*)))
\end{lstlisting}
The input to \cred{subst} is a function $\cblu{\rho} : \N \to \mathsf{\cred{Expr}}$ dictating that \cred{$\oldGamma$Var} $\n$ is to be replaced with $\cblu{\rho}~\n$; we will never substitute for interval variables. When we go under a \cred{Lam} binder, we now want to leave (\cred{$\oldGamma$Var $0$}) untouched, and want to replace (\cred{$\oldGamma$Var $\n$}), for $\n > \cred{0}$, with $\cblu{\rho}~\bkt{\n-\cred{1}}$, except that we also have to weaken this term by adding $\cred{1}$ to all of its \cred{$\oldGamma$Var} variables. Similarly, if we go under a \cred{Line} binder, then we still want to replace \cred{$\oldGamma$Var $\n$} with $\cblu{\rho}~\n$, but we also have to weaken this term by adding $\cred{1}$ to all of its \cred{iVar} variables.

Third, we cover the case of \emph{normalisation}:
\begin{lstlisting}
((*\cred{\textbf{define}}*) ((*\cred{shift}*) (*\cblu{expr}*) (*$\t$*))
  ((*\cred{subst}*) (*\cblu{expr}*) ((*\cred{$\lambda$}*) ((*$\n$*)) ((*\cred{if}*) ((*\cred{zero?}*) (*$\n$*)) (*$\t$*) ((*\cred{$\oldGamma$Var}*) ((*\cred{sub1}*) (*$\n$*)))))))
\end{lstlisting}
\clearpage
\begin{lstlisting}
((*\cred{\textbf{define}}*) ((*\cred{norm}*) (*\cblu{expr}*))
  ((*\cred{\textbf{match}}*) (*\cblu{expr}*)
    (*\cred{$\cdots$}*)
    [((*\cred{App}*) (*$\t$*) (*$\s$*))
     ((*\cred{\textbf{let}}*)
         ([(*$\cblu{tp}$*) ((*\cred{norm}*) (*$\t$*))]
          [(*$\cblu{sp}$*) ((*\cred{norm}*) (*$\s$*))])
       ((*\cred{\textbf{match}}*) (*$\cblu{tp}$*)
         [((*\cred{Lam}*) (*$\cblu{g}$*)) ((*\cred{norm}*) ((*\cred{shift}*) (*$\cblu{g}$*) (*$\cblu{sp}$*)))]
         [(*\cred{else}*) ((*\cred{App}*) (*$\cblu{tp}$*) (*$\cblu{sp}$*))]))]
    (*\cred{$\cdots$}*)))
\end{lstlisting}
We are only including the case of the non-trivial redex. The \cred{shift} helper function substitutes $\t$ in \cblu{expr} for \cred{$\oldGamma$Var $0$}, and shifts down the rest of the variables. To normalise \cred{App} $\t$ $\s$, we normalise $\t$ and $\s$, and check if the normal form of $\t$ is of the form \cred{Lam} $\cblu{g}$. If this is the case, then we call \cred{shift} and recursively normalise the result as new redexes may have been introduced. Otherwise we keep the result as an \cred{App}.

The above has all been standard-fare type theory implementation, but now we get to the dTT specific content. First we cover the \emph{evens renaming operation}:
\begin{lstlisting}
((*\cred{\textbf{define}}*) ((*\cred{evens}*) (*\cblu{expr}*) (*\cblu{$\rho$}*))
  ((*\cred{\textbf{match}}*) (*\cblu{expr}*)
    [((*\cred{XVar}*) (*$\n$*)) ((*\cred{XVar}*) ((*$\cred{\star}$*) (*$\cred{2}$*) (*$\n$*)))]
    [((*\cred{$\mu$Var}*) (*$\n$*)) ((*\cred{$\mu$Var}*) ((*$\cred{\star}$*) (*$\cred{2}$*) (*$\n$*)))]
    [((*\cred{$\oldalpha$Var}*) (*$\n$*)) ((*\cred{$\oldalpha$Var}*) ((*$\cred{\star}$*) (*$\cred{2}$*) (*$\n$*)))]
    [((*\cred{$\oldGamma$Var}*) (*$\n$*)) ((*\cred{$\oldGamma$Var}*) ((*\cblu{$\rho$}*) (*$\n$*)))]
    [((*\cred{iVar}*) (*$\n$*)) ((*\cred{iVar}*) (*$\n$*))]
    [((*\cred{Lam}*) (*$\t$*)) ((*\cred{Lam}*) ((*\cred{evens}*) (*$\t$*) ((*\cred{weaken-ren}*) (*\cblu{$\rho$}*))))]
    [((*\cred{Line}*) (*$\t$*)) ((*\cred{Line}*) ((*\cred{evens}*) (*$\t$*) (*\cblu{$\rho$}*)))]
    (*\cred{$\cdots$}*)))
\end{lstlisting}
This renaming operation behaves uniformly on the $\X$, $\c$, and $\alpha$ variables, but the case for context variables is a bit more difficult. Imagine taking the display of something like $\lambda~\x.~\f~\bkt{\lambda~\y.~\g}$. The result is $\lambda~\x.~\lambda~\x\cblu{'}.~\f\d~\bkt{\bkt{\lambda~\y.~\g}~\sqbkt{\,\x,\x\cblu{'}\,}\ev~}~\bkt{\lambda~\y.~\g}\d$. In this case, the abstraction $\lambda~\y.~\g$ is open. The evens substitution starts off by discarding half of the context variables, but when it recursively goes under the binder on $\y$, it does not discard a variable twinned with $\y$. As a result of this, we keep track of the action on context variables by way of a renaming $\cblu{\rho} : \N \to \N$, which gets weakened when we go under a \cred{Lam} binder.

Now we get to the key function of the implementation -- \emph{display}:
\begin{lstlisting}
((*\cred{\textbf{define}}*) ((*\cred{display}*) (*\cblu{expr}*))
  ((*\cred{\textbf{match}}*)  (*\cblu{expr}*)
    [((*\cred{XVar}*) (*$\n$*)) ((*\cred{XVar}*) ((*\cred{add1}*) ((*$\cred{\star}$*) (*$\cred{2}$*) (*$\n$*))))]
    [((*\cred{$\mu$Var}*) (*$\n$*)) ((*\cred{$\mu$Var}*) ((*\cred{add1}*) ((*$\cred{\star}$*) (*$\cred{2}$*) (*$\n$*))))]
    [((*\cred{$\oldalpha$Var}*) (*$\n$*)) ((*\cred{$\oldalpha$Var}*) ((*\cred{add1}*) ((*$\cred{\star}$*) (*$\cred{2}$*) (*$\n$*))))]
    [((*\cred{$\oldGamma$Var}*) (*$\n$*)) ((*\cred{$\oldGamma$Var}*) ((*$\cred{\star}$*) (*$\cred{2}$*) (*$\n$*)))]
    [((*\cred{Lam}*) (*$\t$*)) ((*\cred{Lam}*) ((*\cred{Lam}*) ((*\cred{display}*) (*$\t$*))))]
    [((*\cred{Line}*) (*$\t$*)) ((*\cred{Line}*) ((*\cred{display}*) (*$\t$*)))]
    [((*\cred{App}*) (*$\t$*) (*$\s$*)) ((*\cred{App}*) ((*\cred{App}*) ((*\cred{display}*) (*$\t$*)) ((*\cred{evens}*) (*$\s$*) ((*\cred{$\lambda$}*) ((*$\n$*)) ((*\cred{add1}*) ((*$\cred{\star}$*) (*$\cred{2}$*) (*$\n$*))))))
                   (*\,*)((*\cred{display}*) (*$\s$*)))]
\end{lstlisting}
\begin{lstlisting}
    [((*\cred{iApp}*) (*$\t$*) (*$\i$*)) ((*\cred{iApp}*) ((*\cred{display}*) (*$\t$*)) (*$\i$*))]
    [((*\cred{Path}*) (*$\X$*) (*$\cblu{e1}$*) (*$\cblu{e2}$*))
     ((*\cred{Lam}*) ((*\cred{PathP}*) ((*\cred{Line}*) ((*\cred{App}*) ((*\cred{rename}*) ((*\cred{display}*) (*$\X$*)) (*\cred{add1}*) (*\cred{add1}*))
                          (*\,*)((*\cred{iApp}*) ((*\cred{$\oldGamma$Var}*) (*$\cred{0}$*)) ((*\cred{iVar}*) (*$\cred{0}$*)))))
                (*\,*)((*\cred{rename}*) ((*\cred{display}*) (*$\cblu{e1}$*)) (*\cred{add1}*) ((*\cred{$\lambda$}*) ((*$\x$*)) (*$\x$*)))
                (*\,*)((*\cred{rename}*) ((*\cred{display}*) (*$\cblu{e2}$*)) (*\cred{add1}*) ((*\cred{$\lambda$}*) ((*$\x$*)) (*$\x$*)))))]
    [((*\cred{PathP}*) ((*\cred{Line}*) (*$\X$*)) (*$\cblu{e1}$*) (*$\cblu{e2}$*))
     ((*\cred{Lam}*) ((*\cred{PathP}*) ((*\cred{Line}*) ((*\cred{App}*) ((*\cred{rename}*) ((*\cred{display}*) (*$\X$*)) (*\cred{add1}*) ((*\cred{$\lambda$}*) ((*$\x$*)) (*$\x$*)))
                          (*\,*)((*\cred{iApp}*) ((*\cred{$\oldGamma$Var}*) (*$\cred{0}$*)) ((*\cred{iVar}*) (*$\cred{0}$*)))))
                (*\,*)((*\cred{rename}*) ((*\cred{display}*) (*$\cblu{e1}$*)) (*\cred{add1}*) ((*\cred{$\lambda$}*) ((*$\x$*)) (*$\x$*)))
                (*\,*)((*\cred{rename}*) ((*\cred{display}*) (*$\cblu{e2}$*)) (*\cred{add1}*) ((*\cred{$\lambda$}*) ((*$\x$*)) (*$\x$*)))))]
    [((*\cred{Comp}*) (*$\cblu{p1}$*) (*$\cblu{p2}$*)) ((*\cred{Comp}*) ((*\cred{display}*) (*$\cblu{p1}$*)) ((*\cred{display}*) (*$\cblu{p2}$*)))]))
\end{lstlisting}
There are two comments to make on this. First, the distinction between de Bruijn levels and indices. For something like \cred{XVar} we index the list of arguments from left to right. If the $\X$ environment consisted of $\sqbkt{\,\X : \Type\,}$, then the displayed environment would look like $\sqbkt{\,\X : \Type,~\X\cblu{'} : \X \to \Type\,}$, and the variable \cred{XVar} $\cred{0}$ should get displayed to \cred{XVar} $\cred{1}$, denoting $\X\cblu{'}$. On other other hand, in the case of the rule for \cred{App}, the first resulting application calls the \cred{evens} renaming operation, and we want this to have the opposite effect by picking the unprimed variable. However, for context variables, we are using indices instead of levels, and thus are indexing from right to left. As a result of this, we use the same $\cred{2}\n+\cred{1}$ formula. Second, the display formulas for $\Path$ and $\PathP$ introduces new binders. Both rules introduce an outermost \cred{Lam} binder, but only the rule for $\Path$ introduces a \emph{new} \cred{Line} binder. Binding new variables means that we must weaken under the binder. In the case of $\X$ in $\Path$ we weaken both the regular and interval contexts, but in the case of $\X$ in $\PathP$, we just weaken the regular context.

From here we can define \emph{meta-abstracted partial substitution d\'ecalage}:
\begin{lstlisting}
((*\cred{\textbf{define}}*) ((*\cred{d\'ecalage}*) (*$\sigma$*))
  ((*\cred{\textbf{match}}*) (*$\sigma$*)
    [(*\cred{`()}*) (*\cred{`()}*)]
    [((*\cred{cons}*) (*$\t$*) (*$\sigma$*)) ((*\cred{cons}*) ((*\cred{evens}*) (*$\t$*) ((*\cred{$\lambda$}*) ((*$\n$*)) ((*\cred{add1}*) ((*$\cred{\star}$*) (*$\cred{2}$*) (*$\n$*)))))
                    (*\,*)((*\cred{cons}*) ((*\cred{display}*) (*$\t$*)) ((*\cred{d\'ecalage}*) (*$\sigma$*))))]))
\end{lstlisting}

Then we write an \emph{evaluation} function that plugs in the meta abstracted arguments when calculating $\cred{\widetilde{\mathsf{S}}}_{\cred{X}}^{\cred{\mathsf{D}}^\k}~\cred{\oldsigma}_{\cred{X}}^\k~\cred{\oldsigma}_{\cred{\mu}}^\k~\cred{\partial a}_\k~\cred{a}_\k$:
\begin{lstlisting}
((*\cred{\textbf{define}}*) ((*\cred{eval}*) (*\cblu{expr}*) (*\cblu{Xenv}*) (*\cblu{$\mu$env}*) (*\cblu{$\oldalpha$env}*))
  ((*\cred{\textbf{match}}*) (*\cblu{expr}*)
    (*\cred{$\cdots$}*)
    [((*\cred{XVar}*) (*$\n$*)) ((*\cred{list-ref}*) (*\cblu{Xenv}*) (*$\n$*))]
    [((*\cred{$\mu$Var}*) (*$\n$*)) ((*\cred{list-ref}*) (*\cblu{$\mu$env}*) (*$\n$*))]
    [((*\cred{$\oldalpha$Var}*) (*$\n$*)) ((*\cred{list-ref}*) (*\cblu{$\oldalpha$env}*) (*$\n$*))]
    (*\cred{$\cdots$}*)))
\end{lstlisting}
The rest of this function is purely structurally recursive and does not modify the \cblu{env} inputs. We are allowed to do so as the environments will consist of closed terms; if they were instead open, we would otherwise have to weaken when going under binders.
\clearpage

Now, we write functions for generating lists of symbols such as $\sqbkt{\,\C,~\A\C,~\B\C,~\A\B\C\,}$:
\begin{lstlisting}
((*\cred{\textbf{define}}*) ((*\cred{update-$\oldalpha$}*) (*$\A$*) (*\cblu{$\oldalpha$gen}*))
  ((*\cred{append}*) (*\cblu{$\oldalpha$gen}*)
         (*\,*)((*\cred{list}*) ((*\cred{list}*) ((*\cred{integer$\to$char}*) ((*$\cred{+}$*) (*$\cred{65}$*) (*$\A$*)))))
         (*\,*)((*\cred{map}*) ((*\cred{$\lambda$}*) ((*$\l$*)) ((*\cred{cons}*) ((*\cred{integer$\to$char}*) ((*$\cred{+}$*) (*$\cred{65}$*) (*$\A$*))) (*$\l$*))) (*\cblu{$\oldalpha$gen}*))))

((*\cred{\textbf{define}}*) ((*\cred{extract-$\oldalpha$}*) (*$\A$*) (*\cblu{$\oldalpha$gen}*))
  ((*\cred{map}*) 
   ((*\cred{$\lambda$}*) ((*$\l$*)) ((*\cred{string$\to$symbol}*) ((*\cred{apply}*) (*\cred{string}*) ((*\cred{reverse}*) (*$\l$*)))))
   ((*\cred{cons}*) ((*\cred{list}*) ((*\cred{integer$\to$char}*) ((*$\cred{+}$*) (*$\cred{65}$*) (*$\A$*))))
        (*\,*)((*\cred{map}*) ((*\cred{$\lambda$}*) ((*$\l$*)) ((*\cred{cons}*) ((*\cred{integer$\to$char}*) ((*$\cred{+}$*) (*$\cred{65}$*) (*$\A$*))) (*$\l$*))) (*\cblu{$\oldalpha$gen}*)))))
\end{lstlisting}

And with all of this, we can finish off our program:
\begin{lstlisting}
((*\cred{\textbf{define}}*) ((*\cred{itr}*) (*$\n$*) (*$\A$*) (*$\cblu{SX}$*) (*$\cblu{S\mu}$*) (*\cblu{Xenv}*) (*\cblu{$\mu$env}*) (*\cblu{$\oldalpha$gen}*))
  ((*\cred{\textbf{cond}}*)
    [((*\cred{zero?}*) (*$\n$*))
     ((*\cred{norm}*)
      ((*\cred{foldl}*) ((*\cred{$\lambda$}*) ((*\cblu{arg}*) (*\cblu{expr}*)) ((*\cred{App}*) (*\cblu{expr}*) (*\cblu{arg}*)))
           (*\,*)((*\cred{last}*) (*\cblu{Xenv}*))
           (*\,*)((*\cred{reverse}*) ((*\cred{rest}*) ((*\cred{reverse}*) ((*\cred{extract-$\oldalpha$}*) (*$\A$*) (*\cblu{$\oldalpha$gen}*)))))))]
    [(*\cred{else}*)
     ((*\cred{\textbf{let}}*) ([(*\cblu{$\oldalpha$batch}*) ((*\cred{extract-$\oldalpha$}*) (*$\A$*) (*\cblu{$\oldalpha$gen}*))])
       ((*\cred{itr}*)
        ((*\cred{sub1}*) (*$\n$*))
        ((*\cred{add1}*) (*$\A$*))
        ((*\cred{d\'ecalage}*) (*$\cblu{SX}$*))
        ((*\cred{d\'ecalage}*) (*$\cblu{S\mu}$*))
        ((*\cred{append}*) (*\cblu{Xenv}*) ((*\cred{map}*) ((*\cred{$\lambda$}*) ((*$\X$*)) ((*\cred{norm}*) ((*\cred{eval}*) (*$\X$*) (*\cblu{Xenv}*) (*\cblu{$\mu$env}*) (*\cblu{$\oldalpha$batch}*)))) (*$\cblu{SX}$*)))
        ((*\cred{append}*) (*\cblu{$\mu$env}*) ((*\cred{map}*) ((*\cred{$\lambda$}*) ((*$\c$*)) ((*\cred{norm}*) ((*\cred{eval}*) (*$\c$*) (*\cblu{Xenv}*) (*\cblu{$\mu$env}*) (*\cblu{$\oldalpha$batch}*)))) (*$\cblu{S\mu}$*)))
        ((*\cred{update-$\oldalpha$}*) (*$\A$*) (*\cblu{$\oldalpha$gen}*))))]))

((*\cred{\textbf{define}}*) ((*\cred{Xn}*) (*$\n$*))
  ((*\cred{itr}*) (*$\n$*)
      (*$\cred{0}$*)
      ((*\cred{list}*) ((*\cred{Lam}*) ((*\cred{Path}*) ((*\cred{XVar}*) (*$\cred{0}$*))
                    (*\,*)((*\cred{App}*) ((*\cred{App}*) ((*\cred{$\mu$Var}*) (*$\cred{0}$*)) ((*\cred{$\oldalpha$Var}*) (*$\cred{0}$*))) ((*\cred{$\oldGamma$Var}*) (*$\cred{0}$*)))
                    (*\,*)((*\cred{App}*) ((*\cred{App}*) ((*\cred{$\mu$Var}*) (*$\cred{0}$*)) ((*\cred{$\oldGamma$Var}*) (*$\cred{0}$*))) ((*\cred{$\oldalpha$Var}*) (*$\cred{0}$*))))))
      ((*\cred{list}*) ((*\cred{Lam}*) ((*\cred{Lam}*) ((*\cred{Lam}*) ((*\cred{Lam}*)
          (*\,*)((*\cred{Comp}*) ((*\cred{Line}*) ((*\cred{App}*) ((*\cred{App}*) ((*\cred{$\mu$Var}*) (*$\cred{0}$*)) ((*\cred{iApp}*) ((*\cred{$\oldGamma$Var}*) (*$\cred{2}$*)) ((*\cred{iVar}*) (*$\cred{0}$*)))) ((*\cred{$\oldGamma$Var}*) (*$\cred{1}$*))))
                (*\;*)((*\cred{Line}*) ((*\cred{App}*) ((*\cred{App}*) ((*\cred{$\mu$Var}*) (*$\cred{0}$*)) ((*\cred{$\oldGamma$Var}*) (*$\cred{3}$*))) ((*\cred{iApp}*) ((*\cred{$\oldGamma$Var}*) (*$\cred{0}$*)) ((*\cred{iVar}*) (*$\cred{0}$*)))))
                (*\;*)))))))
      ((*\cred{list}*) (*$\cred{`X}$*))
      ((*\cred{list}*) (*$\cred{`\mu}$*))
      (*\cred{`()}*)))
\end{lstlisting}
The zero branch of \cred{itr} extracts the top dimensional $\X$ term, iteratively applies it to $\cred{\partial a}_\n$, and normalises. The successor branch implements the formulas for $\cred{\oldsigma}_{\cred{X}}^{\k+\cred{1}}$ and $\cred{\oldsigma}_{\cred{\mu}}^{\k+\cred{1}}$. The defenition of \cred{Xn} instantiates \cred{itr} at $\cred{\widetilde{\S}_X}$ $\cred{\widetilde{\S}_\mu}$, $\cred{\oldsigma_{X}^{0}}$, and $\cred{\oldsigma_{\mu}^{0}}$, thus returning the $\n$-simplex type.

This code, along with a \emph{pretty-print function}, may be found at:
\[\text{\url{https://github.com/FrozenWinters/perm}} \tag*{\bbox}\]
\end{document}